\documentclass[3p,times]{elsarticle}

\usepackage{array}
\usepackage{color}
\usepackage{tabularx}
\usepackage{graphicx}
\usepackage{amsmath}
\usepackage{amssymb}
\usepackage{amsfonts}	
\usepackage{moreverb}
\usepackage{dsfont}
\usepackage{bm}
\usepackage{multirow}
\usepackage{soul}
\usepackage{microtype} 

\newcommand\BibTeX{{\rmfamily B\kern-.05em \textsc{i\kern-.025em b}\kern-.08em
T\kern-.1667em\lower.7ex\hbox{E}\kern-.125emX}}

\newcommand{\x}{\mbf{x}}

\newcommand{\mbf}[1]{\mathbf{#1}}			%

\newcommand{\Q}{\mathbf{Q}}
\renewcommand{\u}{\mathbf{u}}

\newcommand{\q}{\mathbf{q}}
\newcommand{\F}{\mathbf{F}}
\newcommand{\f}{\mathbf{f}}
\newcommand{\g}{\mathbf{g}}
\newcommand{\h}{\mathbf{h}}
\renewcommand{\v}{\mathbf{v}}

\newcommand{\V}{\mathbf{V}}

\newcommand{\bdm}{\begin{displaymath}}
\newcommand{\edm}{\end{displaymath}}

\newcommand{\be}{\begin{equation} \begin{aligned} }
\newcommand{\ee}{\end{aligned} \end{equation}}

\renewcommand{\epsilon}{\varepsilon}
\renewcommand{\phi}{\varphi}
\newcommand{\BoldXi}{\boldsymbol{\xi}}

\newcommand{\bb}[1]{\mbf{#1}}

\newcommand{\dd}{\textup{d}}
\newcommand{\dbl}{\left\llbracket}
\newcommand{\dbr}{\right\rrbracket}

\usepackage[dvipsnames]{xcolor}

\newcommand{\RIIcolor}[1]{{\leavevmode\color{black} #1}}

\newfont{\numerikEleven}{ecrm1000}
\newfont{\numerikTen}{cmss10}
\newfont{\numerikNine}{cmss9}
\newfont{\numerikEight}{cmss8}

\journal{Computers and Fluids}

\begin{document} 
\begin{frontmatter}
\title{A simple diffuse interface approach for compressible flows around moving solids 
	of arbitrary shape based on a reduced Baer-Nunziato model} 

\author[BTU]{Friedemann Kemm}
\ead{friedemann.kemm@b-tu.de}

\author[UniTN]{Elena Gaburro}
\ead{elena.gaburro@unitn.it}

\author[OVG]{Ferdinand Thein}
\ead{ferdinand.thein@ovgu.de}

\author[UniTN]{Michael Dumbser$^{*}$}
\ead{michael.dumbser@unitn.it}
\cortext[cor1]{Corresponding author}

\address[BTU]{BTU Cottbus-Senftenberg, Platz der Deutschen Einheit 1, D-03042 Cottbus, Germany}  
\address[OVG]{Otto-von-Guericke Universit\"at, Universit\"atsplatz 2, D-39106 Magdeburg, Germany}  
\address[UniTN]{Department of Civil, Environmental and Mechanical Engineering, University of Trento, Via Mesiano, 77 - I-38123 Trento, Italy.}

\begin{abstract}
In this paper we propose a new diffuse interface model for the numerical simulation of \textit{inviscid} compressible flows around 
fixed and moving solid bodies of arbitrary shape. The solids are assumed to be moving rigid bodies, without any elastic properties. The mathematical model is a simplified case of the seven-equation  Baer-Nunziato model of compressible multi-phase flows. The resulting governing PDE system is a nonlinear system of hyperbolic conservation laws with non-conservative products. The geometry of the solid bodies is simply specified via a scalar field that represents the volume fraction of the fluid present in each control volume. This allows the discretization of arbitrarily complex geometries on simple uniform or adaptive Cartesian meshes. Inside the solid bodies, the fluid volume fraction is zero, while it is unitary inside the fluid phase. Due to the diffuse interface nature of the model, the volume fraction function can assume any value between zero and one in mixed cells that are occupied by both, fluid and solid. 

We also prove that at the material interface, i.e. where the volume fraction jumps from unity to zero, the normal component of the fluid velocity assumes the value of the normal component of the solid velocity. 
This result can be directly derived from the governing equations, either via Riemann invariants or from the generalized Rankine Hugoniot conditions according to the theory of Dal Maso, Le Floch and Murat \cite{DLMtheory}, which justifies the use of a path-conservative approach for treating the non-conservative products.

The governing partial differential equations of our new model are solved on simple uniform Cartesian grids via a high order path-conservative ADER discontinuous Galerkin (DG) finite element method with \textit{a posteriori} sub-cell finite volume (FV) limiter. Since the numerical method is of the shock capturing type, the fluid-solid boundary is never explicitly tracked by the numerical method, neither via interface reconstruction, nor via mesh motion. 

The effectiveness of the proposed approach is tested on a set of different numerical test problems, including 1D Riemann problems as well as supersonic flows over fixed and moving rigid bodies.  
\end{abstract}

\end{frontmatter}

\section{Introduction} 
\label{sec.intro}
The numerical simulation of fluid-structure-interaction problems in moving compressible media 
is a very important, but at the same time also highly challenging topic. There are overall 
three  different big families of numerical methods in order to tackle this type of problems: 
i) Lagrangian and  Arbitrary-Lagrangian-Eulerian (ALE) methods on moving meshes, where the material interface is exactly resolved and tracked by the moving computational grid; 
ii) Eulerian sharp interface methods on fixed meshes with explicit interface reconstruction, such as the volume of fluid (VOF) method \cite{HirtNichols} or the level-set approach \cite{levelset1,levelset2} in combination with the ghost-fluid method \cite{FedkiwEtAl1,FedkiwEtAl2}; 
iii) Eulerian diffuse interface methods on fixed grids, where the presence of each material is only represented via a scalar color function and where no explicit interface reconstruction technique is applied, see e.g.  \cite{BaerNunziato1986,SaurelAbgrall,Gavrilyuk1,Gavrilyuk2}. 

Probably the most natural choice seems to be a numerical scheme on moving boundary-fitted meshes, where the shape of the solid body is precisely represented by the moving computational grid and conventional wall boundary conditions can be applied at the fluid-solid interface. There is a vast literature on the topic, and it would be impossible to give a complete overview here. Concerning staggered and cell-centered Lagrangian finite volume schemes on moving meshes, we refer the reader to \cite{ShashkovCellCentered,ShashkovRemap3,ShashkovRemap4,ShashkovRemap5,Maire2007,Maire2008,Maire2009,Maire2010,Maire2011,LoubereSedov3D,MaireMM1,MaireMM2,MaireMM3,chengshu1,chengshu2,re2018assessment,ciallella2019shifted,GaburroNonConforming,gaburro2018B} and references therein. Concerning high order purely Lagrangian DG schemes, see the methods forwarded in \cite{Vilar2,Vilar3,Yuetal}.    
In the context of moving mesh schemes, we also mention the family of direct Arbitrary-Lagrangian-Eulerian (ALE) schemes, see for example  \cite{Feistauer1,Feistauer2,Feistauer3,Feistauer4,spacetimedg1,spacetimedg2,LagrangeDG,DGCWENO,GaburroAREPO,gaburroReview} for high order discontinuous Galerkin ALE schemes on moving meshes.  

Concerning an overview of diffuse interface models and related numerical methods, the reader is referred to 
\cite{BaerNunziato1986,SaurelAbgrall,SaurelAbgrall2,SaurelGavrilyukRenaud,SaurelPetitpas,SaurelPetitpasAbgrall,Kapila2001,AndrianovSaurelWarnecke,AndrianovWarnecke,DeledicquePapalexandris,Schwendeman,TokarevaToro,USFORCE2,Menshov2018,Pelanti2006,Pelanti2008,Pelanti2018,Pelanti2019,DIM2D,DIM3D,Gaburro2018,FrontierADERGPR}. 
A diffuse interface model for the interaction of compressible fluids with compressible elasto-plastic solids was 
recently forwarded by the group of Gavrilyuk and Favrie et al. in a series of papers, see \cite{Gavrilyuk1,Gavrilyuk2,Gavrilyuk2008,FavrieGavrilyukSaurel,NdanouFavrieGavrilyuk}.  
For alternative approaches, see also  \cite{Iollo2017,AbateIolloPuppo,Michael2018a,Jackson2019,Barton2019}. 
The models developed and used in the aforementioned references can be considered as complete,
since they fully describe the interaction of compressible flows with compressible elasto-plastic
media. However, there are applications where the elastic deformations of the solid body are not relevant for the computation of the flow field in the fluid, hence these models would result in excessive computational cost and complexity. The objective of the present paper is therefore to derive a \textit{simple} and \textit{reduced} multi-phase flow model based on the diffuse 
interface approach, which is able to describe compressible flows around fixed and moving \textit{rigid solid bodies}.  
An important consequence of this hypothesis will be that the resulting governing PDE system
becomes very simple and easy to solve. Compared to the standard compressible Euler equations, 
there will be only one additional advection equation for the fluid volume fraction, together
with some non-conservative terms that describe the interaction between the fluid and the solid. 

It has to be mentioned that in the context of incompressible flows, another way to embed geometrically complex moving solid obstacles is the so-called immersed boundary method (IBM), 
which goes back to the seminal work of Peskin, see \cite{IBM01}. For an overview of recent developments, see \cite{IBM02,IBM03,IBM04,kim2006immersed,roma1999adaptive,de2006immersed,boukharfane2018combined}
and references therein. 
In immersed boundary methods, an additional force term is added to the momentum equation of the Navier-Stokes equations  which accounts for the presence of the solid body. The non-conservative terms that appear 
in the model derived in this paper will play a similar role as the additional forcing terms in the IBM approach.  An immersed boundary method based on the volume fraction was proposed in~\cite{VolFracIBM}  
and is related to the compressible model presented in this work.  At this point, it is also  important to mention the work by Menshov et al. concerning the description of compressible flows
around complex-shaped objects, see \cite{Menshov2014}. 

The rest of this article is structured as follows: 
in Section~\ref{sec.model}, we detail the derivation of our diffuse interface model describing compressible flows around moving solid obstacles,  
and in particular we provide a detailed proof that the solid and the gas velocities are equal at the material interface; 
then in Section~\ref{sec.method}, we briefly describe the numerical scheme employed for our simulations based on high order ADER-DG schemes with \textit{a posteriori} sub--cell finite volume limiter for dealing with discontinuities. The non-conservative products are treated via the path-conservative approach of Par\'es and Castro \cite{Castro2006,Pares2006,NCproblems, GaburroMNRAS}. 
The obtained numerical results are shown in Section~\ref{sec.results}, and finally, in  Section~\ref{sec.conclusion}, we give some concluding remarks and an outlook to future research and developments. 

\section{Diffuse interface method based on a reduced Baer-Nunziato model} 
\label{sec.model}

To derive the model we start from the full seven equation Baer-Nunziato (BN) model~\cite{BaerNunziato1986,SaurelAbgrall,AndrianovWarnecke,Schwendeman,Kapila2001,MurroneGuillard,Menshov2018} without relaxation source terms, which reads  
\be
\label{eqn.pde.bn} 
&\frac{\partial}{\partial t}\alpha_1+\textbf{v}_I \cdot \nabla \alpha_1 =  0, \\
&\frac{\partial}{\partial t}\left(\alpha_1 \rho_1 \right) 
+ \nabla\cdot\left( \alpha_1 \rho_1 \mathbf{v}_1 \right)  =  0,  \\
&\frac{\partial}{\partial t}\left( \alpha_1 \rho_1 \mathbf{v}_1 \right) 
+\nabla \cdot \left(\alpha_1 \left( \rho_1 \mathbf{v}_1 \otimes \mathbf{v}_1 + p_1 \mathbf{I} \right) \right) 
- p_I \nabla \alpha_1   =  0,  \\ 
&\frac{\partial}{\partial t}\left(\alpha_1 \rho_1 E_1 \right) 
+\nabla \cdot \left[ \alpha_1 \left( \rho_1 E_1 + p_1\right) \mathbf{v}_1 \right] 
- p_I \mathbf{v}_I \cdot \nabla \alpha_1  =  0,  \\
&\frac{\partial}{\partial t}\left(\alpha_2 \rho_2 \right) + \nabla \cdot \left(\alpha_2 \rho_2 \mathbf{v}_2 \right)  = 0,  \\ 
&\frac{\partial}{\partial t}\left( \alpha_2 \rho_2 \mathbf{v}_2 \right) 
+\nabla \cdot \left(\alpha_2 \left( \rho_2 \mathbf{v}_2 \otimes \mathbf{v}_2 + p_2 \mathbf{I} \right) \right) 
- p_I \nabla \alpha_2  =  0,  \\
&\frac{\partial}{\partial t}\left(\alpha_2 \rho_2 E_2 \right) 
+\nabla \cdot \left[ \alpha_2 \left( \rho_2 E_2 + p_2\right) \mathbf{v}_2 \right]  
- p_I \mathbf{v}_I \cdot \nabla \alpha_2  =  0.  \\  
\ee
In the above PDE system $\alpha_j$ denotes the volume fraction of phase number $j$, with $j \in \left\{1,2\right\}$, and the constraint $\alpha_1 + \alpha_2 = 1$. 
Furthermore, $\rho_j$, $\mathbf{v}_j$, $p_j$ and $\rho_j E_j$ represent the density, the velocity vector, the pressure and the total energy per unit mass for 
phase number $j$, respectively. 
Alternatively, the first phase is also called the gas phase (index $g$) and the second phase the solid phase (index $s$), respectively.  

The model \eqref{eqn.pde.bn} is closed by an equation of state (EOS) for each phase $j$ of the form 
\begin{equation}
  e_j = e_j(\rho_j,p_j). 
\end{equation}
The definition of the total energy density 
for each phase is given by 
\be
\rho_j E_j=\rho_j e_j + \frac{1}{2} \rho_j \mathbf{v}_j^2, 
\ee
where $e_i$ is the internal energy.
We further have
\be
\left(\frac{\partial e_i}{\partial\rho_i}\right)_{s_i} = \frac{p_i}{\rho_i^2},\quad \left(\frac{\partial e_i}{\partial s_i}\right)_{\rho_i} = T_i\quad\text{and}\quad
\left(\frac{\partial p_i}{\partial\rho_i}\right)_{s_i} &= \frac{\partial p_i}{\partial\rho_i} = a_i^2, \label{eq:sound_speed_def}
\ee
where $a_i$ denotes the speed of sound of phase $i$ and $s_i$ its specific entropy. For a thermodynamically consistent equation of state these derivatives are well defined and the speed of sound is positive.

For the numerical test problems shown later, we will use the stiffened gas EOS 
\be
e_j=\frac{p_j+\gamma_j \, \pi_k}{\rho_j(\gamma_j-1)}
\ee
with $\gamma_j$ being the ratio of specific heats and $\pi_j$ is a material constant.

\smallskip

In this paper, we choose 
$\textbf{v}_I$ = $\textbf{v}_2$ for the interface velocity and the interface pressure is assumed to be ${p}_{I}$ = ${p}_{1}$. 
This corresponds to the original choice proposed in~\cite{BaerNunziato1986}, which has also been adopted in~\cite{AndrianovWarnecke,Schwendeman,DeledicquePapalexandris,USFORCE2,OsherNC,AMR3DNC}. 
However, alternative choices are also possible, see~\cite{SaurelAbgrall,SaurelAbgrall2}. 

By assuming that the solid phase is the second one and neglecting its elastic deformations, we can therefore consider only rigid body motion of the solid in a given velocity field. With the choice 
$p_I = p_1$ and $\mathbf{v}_I = \mathbf{v}_2$ a \textit{reduced} BN model, similar to the approach presented in~\cite{DIM2D,DIM3D,Gaburro2018,Tavelli2019,FrontierADERGPR} therefore reads 

\be
\label{eqn.pde.red} 
& \frac{\partial}{\partial t}\alpha + \textbf{v}_s \cdot \nabla \alpha  =  0,  \\ 
& \frac{\partial}{\partial t}\left(\alpha \rho \right) 
+ \nabla\cdot\left( \alpha \rho \mathbf{v} \right)  =  0,  \\
& \frac{\partial}{\partial t}\left( \alpha \rho \mathbf{v} \right) 
+\nabla \cdot \left(\alpha \rho \mathbf{v} \otimes \mathbf{v} + \alpha p \, \mathbf{I} \right) 
- p \nabla \alpha   =  0,  \\ 
& \frac{\partial}{\partial t}\left(\alpha \rho E \right) 
+\nabla \cdot \left[ \left( \alpha \rho E + \alpha p \right) \mathbf{v} \right] 
- p \, \mathbf{v}_s \cdot \nabla \alpha  =  0,  \\
& \frac{\partial}{\partial t} \mathbf{v}_s = 0. 
\ee
From now on, for notational simplicity, we will \textit{drop} the subscript $_1$ of the gas phase and only retain the subscript $_s$ for the velocity field of the solid phase. Note that the role of the term $-p \nabla \alpha$ in the momentum equation is similar to the one of the forcing term in immersed boundary methods. In case of a jump in alpha from zero to unity, $\nabla \alpha$ would be the derivative of the Heaviside step function and thus a Dirac delta distribution. However, since we use a \textit{diffuse interface} approach, where the discrete representation of $\alpha$ is usually \textit{smoothed} by numerical dissipation, the term $\nabla \alpha$ will only be an \textit{approximation} of the Dirac distribution. \textcolor{black}{Also note that in regions where $\alpha$ tends from unity to zero the gradient $\nabla \alpha$ in \eqref{eqn.pde.red} naturally plays the role of a \textit{normal vector} to the body surface.  }   

The above system can be written in more compact matrix-vector notation as 
\begin{equation}
   \partial_t \mathbf{Q} + \nabla \cdot \mathbf{F}(\mathbf{Q}) + \mathbf{B}(\mathbf{Q}) \cdot \nabla \mathbf{Q} = \mathbf{0},
   \label{eqn.nc} 
\end{equation}
where $\mathbf{Q} \in \Omega_\mathbf{Q} \subset \mathds{R}^m$ is the vector of conservative variables, $\Omega_\mathbf{Q}$ is the state space, $\mathbf{F} = \mathbf{F}(\mathbf{Q})$ is the nonlinear flux tensor and $\mathbf{B}(\mathbf{Q}) \cdot \nabla \mathbf{Q}$ is a so-called non-conservative 
product. The system \eqref{eqn.nc} is called \textit{hyperbolic} if for all directions $\mathbf{n} \neq \mathbf{0}$ 
the matrix 
\begin{equation*}
\mathbf{A}_n = \left( \partial \mathbf{F} / \partial{\mathbf{Q}} + \mathbf{B} \right) \cdot \mathbf{n}
\label{eqn.A.sys} 
\end{equation*}
has $m$ real eigenvalues and a full set of $m$ linearly independent eigenvectors. 

The proposed model \eqref{eqn.pde.red} allows the representation of moving rigid solid bodies of arbitrarily complex shape on uniform or adaptive Cartesian meshes simply at the aid of the scalar volume fraction function $\alpha$, which is set to $\alpha=0$ inside the solid and to $\alpha=1$ inside the compressible gas. This completely removes the classical mesh generation problem, which can become very cumbersome and time consuming for complex geometries. 

\textit{In the new approach presented in this paper, only Cartesian meshes are used. The entire information related to the geometry of the problem is contained in the scalar function $\alpha$ and is automatically treated via the governing PDE system. \textcolor{black}{No further explicit calculations concerning the geometry of the immersed solid bodies, such as normal vectors, volumes or areas, are needed.}  } 
\textcolor{black}{ This is a unique feature of diffuse interface methods for fluid-structure interaction problems 
and has been used for the first time in the work of Favrie and Gavrilyuk et al.  \cite{Gavrilyuk1,Gavrilyuk2,Gavrilyuk2008,FavrieGavrilyukSaurel,NdanouFavrieGavrilyuk}.
\\
Note that despite the use of high order ADER-DG schemes, the proposed diffuse interface approach on Cartesian grids is necessarily \textit{less accurate} than high order ALE schemes on moving body-fitted meshes, see \cite{Feistauer1,Feistauer2,Feistauer3,Feistauer4,spacetimedg1}. However, the diffuse interface method is 
much easier to implement and can in principle also handle cracks and fragmentation of the solid, see  \cite{Gavrilyuk1,Gavrilyuk2,Gavrilyuk2008,FavrieGavrilyukSaurel,NdanouFavrieGavrilyuk}, while such changes of topology would be much more challenging for traditional moving body-fitted meshes, unless topology changes are explicitly allowed, such as in the methods shown in \cite{PFEM1,PFEM2,PFEM3,PFEM4,PFEM5,PFEM6,Springel,GaburroAREPO}.    
}

\subsection{Solid and gas velocities at the material interface}
\label{sec.proof} 

Now, we want to show that the solid and the gas velocities are equal at the material interface if there is a jump of the gas volume fraction function from $\alpha_R = 1$ in the pure gas to $\alpha_L = 0$ in the solid.  
Since it is easy to prove that the governing equations are rotationally invariant, it is sufficient to focus on the one dimensional case. This simplifies the calculations, by considering only the first component $u$ of vector $\v=(u, v)$ and $u_s$ of the vector $\v_s=(u_s, v_s)$, but the result obtained in the following will be valid in any general normal direction $\mathbf{n} \neq \mathbf{0}$ in the multidimensional case.  

Using the following notation for the conserved variables
\be
\label{eq.ConsVar} 
\Q = (q_1,q_2,q_3,q_4,q_5)^T = (\alpha,\alpha\rho,\alpha\rho u,\alpha \rho E,u_s),
\ee 
we can introduce the conservative flux
\be
\label{eq.flux}
\F(\Q) = \begin{pmatrix} 0\\ q_3\\[8pt] \dfrac{q_3^2}{q_2} + q_1p\\[8pt] \dfrac{q_3}{q_2}\left(q_4 + q_1p\right)\\[8pt] 0\end{pmatrix},
\ee
and the non-conservative part reads
\be
\label{eq.NonConMatrix}
\mbf{B}(\Q) = \begin{pmatrix}
u_s     & 0 & 0 & 0 & 0 \\
0       & 0 & 0 & 0 & 0 \\
-p    & 0 & 0 & 0 & 0 \\
-pu_s & 0 & 0 & 0 & 0 \\
0       & 0 & 0 & 0 & 0 
\end{pmatrix}.
\ee
Thus, in one space dimension the system can be written in the following compact form 
\be
\label{eq.GeneralForm1d}
\partial_t \Q + \partial_x \F(\Q) + \mbf{B}(\Q) \, \partial_x \Q = \mathbf{0}.
\ee
Using the system matrix $\mathbf{A}$ introduced in \eqref{eqn.A.sys}, we can rewrite the system in the following quasilinear form 
\be
\label{eq.GeneralFormNonCons1d}
\partial_t \Q + \mbf{A}(\Q) \, \partial_x \Q = \mathbf{0},
\ee
with
\be
\mbf{A}(\Q) = \begin{pmatrix}
%
q_5 & 0 & 0 & 0 & 0\\
0   & 0 & 1 & 0 & 0\\
q_1\dfrac{\partial p}{\partial q_1}   & -\left(\dfrac{q_3}{q_2}\right)^2 + q_1\dfrac{\partial p}{\partial q_2} & 2\dfrac{q_3}{q_2} + q_1\dfrac{\partial p}{\partial q_3}
& q_1\dfrac{\partial p}{\partial q_4} & 0\\
\left(\dfrac{q_3}{q_2} - q_5\right)p + \dfrac{q_1q_3}{q_2}\dfrac{\partial p}{\partial q_1}  & -\dfrac{q_3(q_4 + q_1p)}{q_2^2} + \dfrac{q_1q_3}{q_2}\dfrac{\partial p}{\partial q_2}
& \dfrac{q_4 + q_1p}{q_2} + \dfrac{q_1q_3}{q_2}\dfrac{\partial p}{\partial q_3}             & \dfrac{q_3}{q_2}\left(1 + q_1\dfrac{\partial p}{\partial q_4}\right)  & 0\\
0   & 0 & 0 & 0 & 0
\end{pmatrix}.\label{jacobian_full}
\ee
Using the primitive variables $\mbf{V} = (\alpha,\rho,u,s,u_s)$ the system can be rewritten as 
\be
\label{eq.GeneralFormNonCons1d_2}
\partial_t \V + \mbf{C}(\V) \, \partial_x \V = \mathbf{0},
\ee
with the new system matrix 
\be
\mbf{C}(\mbf{V}) = \frac{\partial \V}{\partial \Q} \, \mathbf{A}\left( \Q(\V) \right) \, \frac{\partial \Q}{\partial \V} = \begin{pmatrix}
u_s                                 & 0                     & 0         & 0                                                                         & 0\\
\dfrac{\rho}{\alpha}(u - u_s)       & u                     & \rho      & 0                                                                         & 0\\
0                                   & \dfrac{a^2}{\rho}     & u         & \dfrac{1}{\rho}\left(\dfrac{\partial p}{\partial s}\right)_{\rho}         & 0\\
0                                   & 0                     & 0         & u                                                                         & 0\\
0                                   & 0                     & 0         & 0                                                                         & 0
\end{pmatrix}.
\ee
The system matrix has the following eigenvalues 
\begin{equation}
\label{eq.eigenvalues}
 \lambda_0 = u_s, \quad 
 \lambda_{1,1} = u - a, \quad 
 \lambda_{1,2} = u, \quad 
 \lambda_{1,3} = u + a, \quad 
 \lambda_2 = 0,
\end{equation}
which do not explicitly depend on the volume fraction function $\alpha$, as in \cite{Tavelli2019}, hence
the geometric complexity of the solid bodies to be described will not explicitly enter into the CFL stability condition on the time step. 

Note that the submatrix $(\mbf{C})_{ij}$ with $i,j \in \{2,3,4\}$ is the system matrix of the Euler equations of compressible gasdynamics.
According to the calculated eigenvalues we obtain (modulo scaling) the following (right) eigenvectors
\be
\label{eq.eigenvectors}
\mbf{R}_0 &=
\begin{pmatrix} 1 \\ -\dfrac{\rho(u - u_s)^2}{\alpha\left((u - u_s)^2 - a^2\right)} \\ \dfrac{(u - u_s)a^2}{\alpha\left((u - u_s)^2 - a^2\right)} \\ 0 \\ 0\end{pmatrix}\!,\quad
\mbf{R}_{1,1} &= \begin{pmatrix} 0 \\ 1 \\ -\dfrac{a}{\rho} \\ 0 \\ 0\end{pmatrix}\!,\quad
\mbf{R}_{1,2} = \begin{pmatrix} 0 \\ 1\\ 0 \\ a^2\left(\dfrac{\partial p}{\partial s}\right)^{-1}_{\rho}\\ 0\end{pmatrix}\!,\quad
\mbf{R}_{1,3} = \begin{pmatrix} 0 \\ 1 \\ \dfrac{a}{\rho} \\ 0 \\ 0\end{pmatrix}\!,\quad
\mbf{R}_2 &= \begin{pmatrix} 0 \\ 0 \\ 0 \\ 0 \\ 1\end{pmatrix}\!.
\ee
We immediately verify that $\alpha$ may only jump across the wave corresponding to $\mbf{R}_0$ and hence this wave corresponds to the material interface.
It is easy to see that this wave is a contact wave.
The eigenvectors $\mbf{R}_{1,1},\mbf{R}_{1,2}$ and $\mbf{R}_{1,3}$ correspond to the standard eigenvectors of the Euler system. Across these waves $\alpha$ and $u_s$ will not change.
As for the Euler case the fields corresponding to $\mbf{R}_{1,1}$ and $\mbf{R}_{1,3}$ are genuine nonlinear, whereas $\mbf{R}_{1,2}$ is linearly degenerated.
The eigenvector $\mbf{R}_2$ is also a contact, across this wave the solid velocity jumps from $u_{s,L}$ to $u_{s,R}$ given by the initial data.
Written in the original variables $\Q$ the eigenvectors are given by
\be
\label{eq.eigenvectorsCons}
\mbf{R}_0 &=
\begin{pmatrix} 1 \\ -\dfrac{\rho a^2}{(u - u_s)^2 - a^2} \\ -\dfrac{\rho u_s a^2}{(u - u_s)^2 - a^2} \\ 0 \\ 0\end{pmatrix}\!,\quad
\mbf{R}_{1,1} &= \begin{pmatrix} 0 \\ 1 \\ u -\alpha a \\ 0 \\ 0\end{pmatrix}\!,\quad
\mbf{R}_{1,2} = \begin{pmatrix} 0 \\ 1\\ u \\ 0\\ 0\end{pmatrix}\!,\quad
\mbf{R}_{1,3} = \begin{pmatrix} 0 \\ 1 \\ u + \alpha a \\ 0 \\ 0\end{pmatrix}\!,\quad
\mbf{R}_2 &= \begin{pmatrix} 0 \\ 0 \\ 0 \\ 0 \\ 1\end{pmatrix}\!.
\ee

In the following we want to consider the Riemann problem for the following given states (being again $\Q = (\alpha,\alpha\rho,\alpha\rho u,\alpha \rho E,u_s)$)
\be
\Q_L = (0,0,0,0,u_{s,L})^T\quad\text{and}\quad\bb{Q}_R = (1,\rho,\rho u,\rho E, u_{s,R})^T.
\ee
The states left and right of the wave $\mathbf{R}_0$ will be denoted with $\mbf{Q}^-$ and $\mbf{Q}^+$, respectively.

\subsubsection{Proof based on Riemann invariants}

After rescaling the eigenvector in~\eqref{eq.eigenvectorsCons}, we have the following useful Riemann invariants for the contact~$\mbf{R}_0$
\be
\frac{\dd q_2}{\dd \sigma} = 1,\quad \frac{\dd q_3}{\dd \sigma} = q_5\quad\text{and}\quad\frac{\dd q_5}{\dd \sigma} = 0,
\ee
\RIIcolor{where $\sigma$ is the independent variable of the parametrization of the corresponding integral curve, see \cite{dafermos2005hyperbolic};} they
can be reformulated to
\be
\label{eq.RiemannInvariantq5}
\frac{\dd q_3}{\dd q_2} = q_5,
\ee
and since {$q_5 \equiv u_s$} is constant across $\mbf{R}_0$, we have
{\be
	\int_{\mbf{Q}^-}^{\mbf{Q}^+}\,\dd q_3 = u_s\int_{\mbf{Q}^-}^{\mbf{Q}^+}\,\dd q_2 .  
	\ee
Due to the structure of the eigenvectors we have $q_1^+ = \alpha_R = 1$ (gas) and $q_1^- = \alpha_L = 0$ (solid) and it follows that 
\be
q_3^+ = q_5q_2^+\quad\Leftrightarrow\quad\rho^+u^+ = u_s\,\rho^+\quad\Leftrightarrow\quad u^+ = u_s.
\ee
}
Thus, we have for Riemann initial data with $\alpha_{L} = 0$ (solid) and $\alpha_{R} = 1$ (gas) that the velocity of the gas $u^+$ at the material interface is equal to the solid velocity $u_s$.
{\scriptsize{$\blacksquare$}}

\subsubsection{Proof based on the generalized Rankine Hugoniot conditions}

As discussed before (below Equation~\eqref{eq.eigenvectors}), we have only one wave where $\alpha$ may change and across the others it remains constant.
In other words the non-conservative products vanish for $\mbf{R}_{1,1},\mbf{R}_{1,2},\mbf{R}_{1,3}$ and $\mbf{R}_2$.
Thus we have a conservative system there and may use the standard Rankine Hugoniot conditions across discontinuities.
For the remaining wave $\mbf{R}_0$ we follow the approach introduced by Dal Maso, Le Floch and Murat in~\cite{DLMtheory}.
Hence, we can establish jump conditions using the formula  
\be
\label{eq.GRH}
S\dbl\bb{Q}\dbr = \int_0^1\bb{A}(\Psi(\bb{Q}^-,\bb{Q}^+,\tau))\frac{\partial\Psi}{\partial\tau}(\tau;\bb{Q}^-,\bb{Q}^+)\,\dd\tau.
\ee
Here $\bb{Q}^-$ and $\bb{Q}^+$ denote the left and right states of the wave $\mbf{R}_0$.
The matrix $\mathbf{A}$ is given by~\eqref{jacobian_full} and $\Psi(\tau;\bb{Q}^-,\bb{Q}^+)$ denotes a suitable path connecting both states properly, see~\cite{DLMtheory}.
In particular we choose the straight-line segment path
\be
\Psi(\bb{Q}^-,\bb{Q}^+,\tau) = (1 - \tau)\bb{Q}^- + \tau\bb{Q}^+.
\ee
We now want to calculate the right state $\bb{Q}^+$ for a given left state $\bb{Q}^-$.
Since $\alpha$ only changes across $\bb{R}_0$, we conclude from the initial data, that $q_1^- = 0$ and $q_1^+ = 1$, 
and we have the following left state
\be
\label{eq.LeftState}
\bb{Q}^- = (0,0,0,0,q_5^-)^T.
\ee
Hence for each $\tau \in [0,1]$ the path defines a state
\be
\bb{Q}_\tau = \begin{pmatrix}\tau q_1^+\\ \tau q_2^+\\ \tau q_3^+\\ \tau q_4^+\\ (1 - \tau)q_5^- + \tau q_5^+\end{pmatrix}
\equiv \begin{pmatrix}\tau\\ \tau\rho^+\\ \tau\rho^+ u^+\\ \tau\rho^+E^+\\ (1 - \tau)u_s^- + \tau u_s^+\end{pmatrix}.
\ee
Further we have
\be
\frac{\partial\Psi}{\partial\tau}(\bb{Q}^-,\bb{Q}^+,\tau) = \Delta \bb{Q} = \begin{pmatrix} q_1^+\\ q_2^+\\ q_3^+\\ q_4^+\\ \Delta q_5\end{pmatrix}
\equiv \begin{pmatrix} 1\\ \rho^+\\ \rho^+ u^+\\ \rho^+E^+\\ \Delta u_s\end{pmatrix},
\ee
\RIIcolor{(where $\Delta \Q = \Q^+ - \Q^-$).} Assuming the pressure to be given as $p = p(\rho,e)$ we have
\be
& e(\bb{Q}) = \frac{q_4}{q_2} - \frac{1}{2}\left(\frac{q_3}{q_2}\right)^2
\quad\text{with}\quad
\dd e = -\left(\frac{q_4}{q_2^2} - \frac{q_3^2}{q_2^3}\right)\,\dd q_2 - \frac{q_3}{q_2^2}\,\dd q_3 + \frac{1}{q_2}\,\dd q_4\\
& \dd p =  -\frac{q_2}{q_1^2}\frac{\partial p}{\partial \rho}\,\dd q_1 + \left(\frac{\partial p}{\partial \rho} + \frac{\partial p}{\partial e}\frac{\partial e}{\partial q_2}\right)\,\dd q_2
+ \frac{\partial p}{\partial e}\frac{\partial e}{\partial q_3}\,\dd q_3 + \frac{\partial p}{\partial e}\frac{\partial e}{\partial q_4}\,\dd q_4.
\ee
Thus we verify that everywhere the $\tau$ cancels except for the terms where $q_5$ occurs, and $\bb{A}(\bb{Q}_\tau)$ literally is the same as~\eqref{jacobian_full} if we replace $q_5$ with
{$u_{s,\tau} = (1 - \tau)u_s^- + \tau u_s^+$} and use the right $+$ values for $q_1,\dots,q_4$.
{Integrating these terms gives $\bar{u}_s = \dfrac{1}{2}(u_s^+ + u_s^-)$.}
Then the integration of~\eqref{eq.GRH} results in (neglecting the $+$ in matrix $\mbf{A}$ in the first passage)
\be
\int_0^1&\bb{A}(\Psi(\bb{Q}^-,\bb{Q}^+,\tau))\frac{\partial\Psi}{\partial\tau}(\bb{Q}^-,\bb{Q}^+,\tau)\,\dd\tau  = \\
= &\begin{pmatrix}
%
\bar{u}_s & 0 & 0 & 0 & 0\\
0   & 0 & 1 & 0 & 0\\
q_1\dfrac{\partial p}{\partial q_1}   & -\left(\dfrac{q_3}{q_2}\right)^2 + q_1\dfrac{\partial p}{\partial q_2} & 2\dfrac{q_3}{q_2} + q_1\dfrac{\partial p}{\partial q_3}
& q_1\dfrac{\partial p}{\partial q_4} & 0\\
\left(\dfrac{q_3}{q_2} - \bar{u}_s\right)p + \dfrac{q_1q_3}{q_2}\dfrac{\partial p}{\partial q_1}  & -\dfrac{q_3(q_4 + q_1p)}{q_2^2} + \dfrac{q_1q_3}{q_2}\dfrac{\partial p}{\partial q_2}
& \dfrac{q_4 + q_1p}{q_2} + \dfrac{q_1q_3}{q_2}\dfrac{\partial p}{\partial q_3}             & \dfrac{q_3}{q_2}\left(1 + q_1\dfrac{\partial p}{\partial q_4}\right)  & 0\\
0   & 0 & 0 & 0 & 0
\end{pmatrix}\cdot\begin{pmatrix} q_1^+\\ q_2^+\\ q_3^+\\ q_4^+\\ \Delta q_5\end{pmatrix}\\
= &\begin{pmatrix} \bar{u}_s q_1^+\\ q_3^+\\
q^+_1\left(q_1^+\dfrac{\partial p}{\partial q_1} + q_2^+\dfrac{\partial p}{\partial q_2} + q_3^+\dfrac{\partial p}{\partial q_3} + q_4^+\dfrac{\partial p}{\partial q_4}\right)
-\dfrac{(q^+_3)^2}{q_2^+} + 2\dfrac{(q^+_3)^2}{q^+_2}\\
q_1^+\left(\dfrac{q^+_3}{q^+_2} - \bar{u}_s \right)p - \dfrac{q^+_3(q^+_4 + q^+_1 p)}{q^+_2} + \dfrac{q_3^+(q^+_4 + q^+_1 p)}{q^+_2} + \dfrac{q^+_3q_4^+}{q^+_2}
+ \dfrac{q^+_1q^+_3}{q^+_2}\left(q_1^+\dfrac{\partial p}{\partial q_1} + q_2^+\dfrac{\partial p}{\partial q_2} + q_3^+\dfrac{\partial p}{\partial q_3} + q_4^+\dfrac{\partial p}{\partial q_4}\right)\\ 0\end{pmatrix}\\
= &\begin{pmatrix} \bar{u}_s q_1^+\\ q_3^+\\ \dfrac{(q_3^+)^2}{q_2^+} - (1 - q_1^+)\left(\dfrac{\partial p}{\partial \rho_1}\right)_{e_1}\\
q_1^+\left(\dfrac{q^+_3}{q^+_2} - \bar{u}_s \right)p + \dfrac{q^+_3q_4^+}{q^+_2} - \dfrac{(1 - q_1^+)q_3^+}{q_2^+}\left(\dfrac{\partial p}{\partial \rho_1}\right)_{e_1}\\ 0\end{pmatrix}.
\ee
Thus $\bb{Q}^+$ is defined by the equation
\be
\label{eq.Qplus1}
\begin{pmatrix} \bar{u}_s q_1^+\\ q_3^+\\ \dfrac{(q_3^+)^2}{q_2^+} - (1 - q_1^+)\left(\dfrac{\partial p}{\partial \rho}\right)_{e}\\
q_1^+\left(\dfrac{q^+_3}{q^+_2} - \bar{u}_s \right)p + \dfrac{q^+_3q_4^+}{q^+_2} - \dfrac{(1 - q_1^+)q_3^+}{q_2^+}\left(\dfrac{\partial p}{\partial \rho}\right)_{e}\\ 0\end{pmatrix}
= S\dbl\bb{Q}\dbr = S\begin{pmatrix} q_1^+\\ q_2^+\\ q_3^+\\ q_4^+\\ \Delta q_5\end{pmatrix}.
\ee
Since we already know from the previous analysis (above Equation~\eqref{eq.LeftState}) that $q_1^+ = 1$,~\eqref{eq.Qplus1} reduces to
\be
\begin{pmatrix} \bar{u}_s \\ q_3^+\\ \dfrac{(q_3^+)^2}{q_2^+}\\ \dfrac{q_3^+q_4^+}{q_2^+} + \left(\dfrac{q_3^+}{q_2^+} - \bar{u}_s \right)p \\ 0\end{pmatrix}
= S\begin{pmatrix} 1\\ \rho^+\\ \rho^+ u^+\\ \rho^+E^+\\ \Delta u_s\end{pmatrix}\!.
\ee
Assuming $S \neq 0$ gives $\Delta u_s = 0$ and thus (since $\bar{u}_s = 1\cdot S$) $S = \bar{u}_s = u_s^- = u_s^+ := u_s$ as expected (indeed it is a contact). 
Further we have $q_3^+ = S q_2^+$ which now gives the desired result
\be 
u_s = \frac{q_3^+}{q_2^+} = u^+.
\ee 
This is the same as the obtained Riemann invariant \eqref{eq.RiemannInvariantq5}.
Finally, for $S = 0$ it is trivial. 
Since we have a contact it follows that this is only the case for $S = u_s^- = u_s^+ = 0$ and again $u^+ = \frac{q_3^+}{q_2^+} = 0 =u_s$.
{\scriptsize{$\blacksquare$}}

\section{Brief summary of the high order path-conservative DG scheme with a posteriori sub-cell finite volume limiter}
\label{sec.method}

As already discussed before, the reduced Baer-Nunziato model~\eqref{eqn.pde.red} in $d$-space dimensions, under consideration in this work, can be cast in the following general form
\be
\frac{\partial \mathbf{Q} }{\partial t} + \nabla \cdot \mathbf{F}\left( \mathbf{Q} \right) 
+ \mathbf{B}(\Q) \cdot \nabla \mathbf{Q}  = \mbf{0}, \qquad \x \in \Omega \subset \mathds{R}^d, \quad t \in \mathds{R}_0^+,
\label{eqn.pde}
\ee
which describes nonlinear systems of hyperbolic equations with non-conservative products, 
and  for which we recall that $\Q(\x, t) \in \Omega_Q \subset \mathds{R}^\nu$ is the state vector of $\nu$ conserved quantities, 
$\F(\Q) = (\f,\g,\h)$ is a non-linear flux tensor, 
$\mathbf{B}(\Q) \cdot \nabla \mathbf{Q} $ collects the non-conservative products,
and $\Omega$ denotes the computational domain,
whereas $\Omega_Q$ is the space of physically admissible states. 
Without loss of generality, we will present the method for the case $d=3$.

To solve~\eqref{eqn.pde} we employ a discontinuous Galerkin (DG) scheme of arbitrary high order of accuracy both in space and in time, based on a fully discrete ADER predictor-corrector procedure, first proposed  in~\cite{DumbserEnauxToro,Dumbser2008} and then detailed for the Cartesian case in~\cite{AMR3DCL,ADERDGVisc,ADERGRMHD,FrontierADERGPR,DumbserGLM,fambri2020discontinuous}. Note that the original ADER approach was introduced by Toro and Titarev in the context of finite volume schemes, using an approximate solution 
of the generalized Riemann problem \cite{menshov} with piecewise polynomial initial data, see \cite{toro3,toro4,titarevtoro,Toro:2006a} for details.  
Here below we briefly summarize the key ingredients of the scheme: after having introduced in Section~\ref{ssec.Domain} the domain discretization and the polynomial data representation $\u_h$, we explain how to evolve $\u_h$ \textit{in the small} (see~\cite{harten}), i.e. without needing of any communications with the neighbors, in order to obtain a \textit{predictor} $\q_h$ of the solution of high order in space and also in \textit{time}. 
This predictor will be then used in the final \textit{corrector} step described in Section~\ref{ssec.corrector}, where the weak form of~\eqref{eqn.pde} is integrated in space and time, and where the numerical fluxes and non-conservative products are evaluated making use of the predictor $\q_h$. The corrector step evolves the discrete solution $\u_h$ in time and takes into account the information coming from the cell neighbors via classical numerical flux functions (approximate Riemann solvers).
We close the section by describing our \textit{a posteriori} sub-cell finite volume limiter~\cite{Dumbser2014,DGLimiter2,DGLimiter3,ADERDGVisc,DGCWENO}, 
which assures the robustness of the scheme even in the presence of discontinuities, but keeping at the same time also the high resolution of the underlying DG scheme.

\subsection{Domain discretization and high order data representation in space}
\label{ssec.Domain}

We discretize $\Omega$ by covering it with a Cartesian grid, called \textit{main grid}, made of $N_E =N_x \times N_y\times N_z$ conforming elements (quadrilaterals if $d=2$, or hexahedra if $d=3$) $\Omega_{ijk}, i = 1, \dots, N_x, j=1, \dots N_y, k=1, \dots N_z$, with volume $|\Omega_{ijk}| = \int_{\Omega_{ijk}} d\x$ and such that
\be
& \Omega_{ijk} \!=\! [x_{i-\frac{1}{2}},x_{i+\frac{1}{2}}]\times [y_{j-\frac{1}{2}},y_{j+\frac{1}{2}}]\times [z_{k-\frac{1}{2}},z_{k+\frac{1}{2}}], \\  
&\text{with } \ 
\Delta x_{i} \!=\! x_{i+\frac{1}{2}}-x_{i-\frac{1}{2}}, \quad 
\Delta y_{j} \!=\! y_{j+\frac{1}{2}}-y_{j-\frac{1}{2}}, \quad 
\Delta z_{k} \!=\! z_{k+\frac{1}{2}}-z_{k-\frac{1}{2}}.
\label{eq.MainGrid}
\ee 
Moreover, for each element we define a reference frame of coordinates $\BoldXi=(\xi, \eta, \zeta)$ linked to the Cartesian coordinates $\x=(x,y,z)$ of $\Omega_{ijk}$ by
\be
\label{eq.mapping}
x = x_{i-\frac{1}{2}} + \xi \Delta x, \quad  y = y_{j-\frac{1}{2}} + \eta \Delta y, \quad
z = z_{k-\frac{1}{2}} + \zeta \Delta z, \quad \xi, \eta, \zeta \in [0,1].
\ee 

Then, in each cell $\Omega_{ijk}$, at the beginning of each time step, the conserved variables $\Q$ are represented at the aid of $d-$dimensional piecewise polynomials of degree $N$ 
\be
\mathbf{u}_h(\x,t^n) = \mathbf{u}_h(\BoldXi(\x)) = \sum \limits_{\ell=0}^{\mathcal{N}-1} \phi_\ell(\BoldXi) \, \hat{\mathbf{u}}_{\ell} 
:= \phi_\ell(\BoldXi) \, \hat{\mathbf{u}}_{\ell} , \quad \x \in \Omega_{ijk}, \quad \mathcal{N} = (N+1)^d,
\label{eqn.uh}
\ee
where $\phi_\ell(\BoldXi)$ are \textit{nodal} spatial basis functions 
given by the tensor product of a set of Lagrange interpolation polynomials of maximum degree $N$ such that
\be
\label{eq.nodalGeneralBasis}
\phi_\ell(\BoldXi_{\text{GL}}^m) = \phi_{\ell_1}\!\left(\xi_{\text{GL}}^m\right)\phi_{\ell_2}\!\left(\eta_{\text{GL}}^m\right)\phi_{\ell_3}\!\left(\zeta_{\text{GL}}^m\right) 
= \left\{ \begin{array}{rl} 1 & \text{if }
	\ell_i=m \\0 & \text{otherwise} \end{array}\right.\hspace{0.7cm}
\ell_i,m=1,\ldots,(N+1),
\ee 
where $\BoldXi_{\text{GL}}^m$ are the set of $(N+1)^d$ Gauss-Legendre (GL)
quadrature points obtained by the tensor product of the GL quadrature points $\xi_{\text{GL}}^m, \eta_{\text{GL}}^m, \zeta_{\text{GL}}^m$ in the unit interval $[0,1]$.

Let us finally underline that the use of a Cartesian grid makes it possible to work in a dimension by dimension fashion, which remarkably reduces the computational cost of the entire algorithm.

\subsection{High order in time via an element-local space-time discontinuous Galerkin predictor}
\label{ssec.predictor}

Representing the conserved variables through high order piecewise polynomials~\eqref{eqn.uh} 
already provides by construction high order of accuracy from the \textit{spatial} point of view. 
Now, in order to achieve also high order of accuracy in \textit{time}, 
we rely on the ADER predictor-corrector approach, 
which strongly differs from conventional semi-discrete Runge-Kutta based methods and which leads to 
a fully-discrete one-step method.  
The predictor step is fully local and avoids any interactions with the neighbors, and is thus well suited
for parallel computing. It also results to be much simpler with respect to the cumbersome Cauchy-Kovalevskaya procedure used in traditional ADER schemes \cite{toro2,toro3,toro4,titarevtoro,Toro:2006a}. 

The so-called predictor $\q_h$ is a space-time polynomial of degree $N$ in $(d+1)$-dimensions which takes the following form 
\be
\q_h(\x, t) = \q_h(\BoldXi(\x), \tau(t)) = \sum_{\ell=0}^{\mathcal{Q}-1} \theta_\ell (\BoldXi, \tau) \hat{\q}_\ell
= \theta_\ell (\BoldXi, \tau) \hat{\q}_\ell,\quad \x \in \Omega_{ijk},\quad t \in [t^n, t^{n+1}],\quad \mathcal{Q} = (N+1)^{d+1},
\label{eqn.qh}
\ee 
where again $\theta_\ell (\BoldXi, \tau)$ is given by the tensor product of Lagrange interpolation polynomials $\phi_{\ell}\left(\BoldXi (\x) \right) \phi_{\ell_\tau}\left(\tau\right)$, with 
$\BoldXi(\x)$ given by~\eqref{eq.mapping} and the mapping for the time coordinate given by $t = t^n + \tau \Delta t, \tau \in [0,1]$.

In order to determine the unknown coefficients $\hat{\q}_\ell$ of~\eqref{eqn.qh} we search $\q_h$ such that it satisfies 
a weak form of the governing PDE~\eqref{eqn.pde} integrated in space and time locally \textit{inside} each $\Omega_{ijk}$ \RIIcolor{(with $\Omega_{ijk}^{\circ} = \Omega_{ijk} \backslash \partial \Omega_{ijk} $ being the interior of $\Omega_{ijk}$)}
\be
\int_{t^n}^{t^{n+1}} \!\! \int_{\Omega_{ijk}^{\circ} } \theta_k \, \partial_t \q_h \,d\mbf{x} \, dt 
+ \int_{t^n}^{t^{n+1}} \!\! \int_{\Omega_{ijk}^{\circ} } \theta_k \,\nabla \cdot \F(\q_h) \,d\x\,dt 
+ \int_{t^n}^{t^{n+1}} \!\! \int_{\Omega_{ijk}^{\circ} } \theta_k \mbf{{B}}(\q_h ) \cdot \nabla \q_h \,d\x\,dt 
= \mbf{0},
\label{eq.predictorEq}
\ee
where the first term is integrated in time by parts exploiting the \emph{causality principle}
(upwinding in time) 
\be
& \int_{\Omega_{ijk}^{\circ}}  \theta_k(\x,t^{n+1}) \q_h(\x,t^{n+1}) \, d\x -
\int_{\Omega_{ijk}^{\circ}} \theta_k(\x,t^{n}) \mbf{u}_h(\x,t^{n}) \, d\x - 
\int_{t^n}^{t^{n+1}} \!\!  \int_{\Omega_{ijk}^{\circ} } \!\!\! \partial_t \theta_k(\x,t) \q_h(\x,t) \,d\x\, dt \\
& + \int_{t^n}^{t^{n+1}} \!\!  \int_{\Omega_{ijk}^{\circ} } \!\!\! \theta_k(\x,t) \nabla \cdot \F(\q_h(\x,t)) \,d\x\,dt 
+ \int_{t^n}^{t^{n+1}} \!\!  \int_{\Omega_{ijk}^{\circ} } \theta_k(\x,t) \mbf{{B}}(\q_h(\x,t)) \cdot \nabla\q_h(\x,t) \,d\x\,dt = \mbf{0},
\label{eq:DOFpredictor}
\ee
and $\mbf{u}_h(\x,t^{n})$ is the known initial condition at time $t^n$. 

Now, the system~\eqref{eq:DOFpredictor}, which contains only volume integrals to be calculated inside $\Omega_{ijk}$ and no surface integrals,  
can be solved via a simple discrete Picard iteration for each element $\Omega_{ijk}$, 
and there is no need of any communication with neighbor elements. 
We recall that this procedure has been introduced for the first time in~\cite{Dumbser2008} for unstructured meshes, 
it has been extended for example to moving meshes in~\cite{Lagrange2D} and
to degenerate space time elements in~\cite{GaburroAREPO}; 
finally, its convergence has been formally proved in~\cite{FrontierADERGPR}.

\subsection{Fully discrete one-step path-conservative ADER-DG scheme}
\label{ssec.corrector}

The update formula of our ADER-DG scheme is recovered, as usual, starting from 
the weak formulation of the governing equations~\eqref{eqn.pde} 
(where the test functions $\phi_k$ coincide with the basis functions $\phi_\ell$ of~\eqref{eq.nodalGeneralBasis})
\be
\int_{t^n}^{t^{n+1}} \int_{\Omega_{ijk}} \phi_k \left(\partial_t \Q  + \nabla \cdot \F(\Q) 
+ \mbf{{B}}(\Q) \cdot\nabla \Q \right) \,d\x\,dt = \mbf{0};
\label{eqn.pde.weak}
\ee
we then substitute $\Q$ with~\eqref{eqn.uh} at time $t=t^n$ (the known initial condition) 
and at $t=t^{n+1}$ (to represent the unknown evolved conserved variables),
and with the high order predictor $\q_h$ previously computed for $t\in[t^n, t^{n+1}]$, obtaining
\be
& \left( \, \int_{\Omega_{ijk}} \phi_k \phi_l \, d\x\right)
\left( \hat{\u}_\ell^{n+1} - \hat{\u}_\ell^{n} \, \right) 
+ \int_{t^n}^{t^{n+1}} \!\! \int_{\partial \Omega_{ijk}} \!\! \phi_k \mathcal{D}\left(\q_h^-, \q_h^+ \right) \cdot \mbf{n} \, dS \, dt \\
& - \int_{t^n}^{t^{n+1}} \!\! \int_{\Omega_{ijk} } \!\!\! \nabla \phi_k \cdot \F(\q_h) \,d\x\,dt 
+ \int_{t^n}^{t^{n+1}} \!\! \int_{\Omega_{ijk}^{\circ} } \phi_k \mbf{{B}}(\q_h) \cdot \nabla \q_h \,d\x\,dt = \mbf{0}.
\label{eq:ADER-DG}
\ee
The use of $\q_h$ allows to compute the integrals appearing in~\eqref{eq:ADER-DG} with high order of accuracy.
Moreover, note that due to the discontinuous character of the solution $\q_h$ at the element interfaces $\partial \Omega_{ijk}$,
the jump term $\mathcal{D}$, which contains the numerical flux as well as a discretization of the non-conservative product, is computed through a numerical flux function 
evaluated over $\q_h^{-}$ and $\q_h^{+}$ which are the so-called boundary-extrapolated data.
In particular, we have employed a two-point path-conservative numerical flux function of Rusanov-type which reads as follows 
\be 
\label{eqn.fluxPC}
\mathcal{D}(\q_h^{-},\q_h^{+}) \cdot \mathbf{n} 
= \ \frac{1}{2} \left( {\mbf{F}}(\q_h^{+}) + {\mbf{F}}(\q_h^{-}) \right) \cdot {\mbf{n}} 
- \frac{1}{2} s_{\max} \left( \q_h^{+} - \q_h^{-} \right) + \frac{1}{2} \left(\int \limits_{0}^{1} {\mbf{B}} \left(\mbf{\Psi}(\q_h^{-},\q_h^{+},\tau) \right)\cdot\mbf{n} \, ds \right)
\cdot\left(\q_h^{+} - \q_h^{-}\right),
\ee
where $s_{\max}$ is the maximum eigenvalue of the system matrices $\mathbf{A}(\q_h^{+})$ and $\mathbf{A}(\q_h^{-})$ 
being
\be 
\label{eq.ALEjacobianMatrix}
\mathbf{A}(\Q)=\frac{\partial \mathbf{F}}{\partial \Q} + \mbf{B},
\ee
and the path $\mbf{\Psi}=\mbf{\Psi}(\q_h^-,\q_h^+, s)$ is the straight-line segment path
\be
\mbf{\psi} = \mbf{\psi}(\q_h^-, \q_h^+, s) = \q_h^- + s \left( \q_h^+ - \q_h^- \right)\,,  \qquad  s \in [0,1]\,, 
\ee
connecting $\q_h^{-}$ and $\q_h^{+}$ in phase-space. The path allows to treat the jump of the non-conservative 
products according to the theory introduced by Dal Maso, Le Floch and Murat in \cite{DLMtheory} (DLM theory) 
and is used for the construction of so-called path-conservative schemes, see  \cite{Pares2004,Pares2006,Castro2006,Munoz2007,Castro2d,NCproblems} for details. For the extension of path-conservative schemes to DG and finite volume methods of arbitrary high order, see \cite{Rhebergen2008,ADERNC}. As already shown before, the straight line segment path is also consistent with the boundary condition that requires the local fluid velocity to be equal to the local solid velocity when the volume fraction jumps from unity to zero. 
\textcolor{black}{Note that the segment path is merely needed for the definition of the jump terms at the element boundaries in the presence of non-conservative products of the type $\mathbf{B}(\mathbf{Q}) \cdot \nabla \mathbf{Q}$ in the framework of path-conservative schemes and according to the DLM theory \cite{DLMtheory}; it has nothing to do with a piecewise linear representation of the geometry of the solid bodies. The geometry of the solids is only represented by the spatial distribution of the scalar field $\alpha$. } 

We conclude this section by recalling that, since the employed method is explicit, the time step $\Delta t$ 
is computed under a classical (global) Courant-Friedrichs-Levy (CFL) stability condition with 
$\textnormal{CFL} \leq 1$ and it is given by 
\be 
\Delta t_\text{DG} < \text{CFL}\frac{h_{\text{min}} }{d\left(2N+1\right)} \frac{1}{|\lambda_{\text{max}}|}
\label{eq.timestep}
\ee
where $h_{\text{min}}$ is the minimum characteristic mesh-size and $|\lambda_{\max}|$ is the spectral radius of the system matrix $\mathbf{A}$.

\subsubsection{Adaptive mesh refinement (AMR)}

Furthermore, in order to increase the resolution in the areas of interest, 
the ADER-DG scheme described above has been implemented on space-time adaptive Cartesian meshes, 
with a \textit{cell-by-cell} refinement approach; for all the details we refer to~\cite{AMR3DCL,AMR3DNC,DGLimiter2,ADERDGVisc,ADERGRMHD,Peano1,Peano2,Exahype}.

The main idea behind our AMR technique consists in, starting from the main grid~\eqref{eq.MainGrid}, 
introducing successive refinement levels, built according to the so
called refinement factor $\mathfrak{r}$, which is the number of sub elements 
per space-direction in which a coarser element is broken
in a refinement process, or which are merged in a recoarsening stage. 
The refinement/recoarsening process is driven by a prescribed {refinement-estimator function} well described in the above references.
Finally, the numerical solution at the sub-cell level during a refinement step is obtained by a standard $L_2$ projection, 
while a reconstruction operator is employed to recover the solution on the main grid starting from the sub-cell level.
Projection~\eqref{eq.projection} and reconstruction~\eqref{eq.reconstruction}-\eqref{eqn.LSQ} are also used in the limiter procedure and hence are better described in next Section~\ref{ssec.limiter}.

\subsection{\textit{A posteriori} sub-cell finite volume limiter}
\label{ssec.limiter}

Higher order discontinous Galerkin schemes can be seen as linear schemes in the sense of Godunov \cite{godunov}, 
hence, in presence of discontinuities, spurious oscillations typically arise.
To minimize their effects we adopt an \textit{a~posteriori} limiting procedure based on the MOOD paradigm 
\cite{MOOD,MOODhighorder,ADERMOOD}: indeed, we first apply our unlimited ADER-DG scheme everywhere, and then, at the end of each time step, we check \textit{a posteriori} the reliability of the obtained candidate solution $\u_h^*$ in each cell.
This candidate solution is checked against physical and numerical admissibility criteria, 
such as floating point exceptions, violation of positivity or other physical bounds, or violation of a relaxed discrete maximum principle (DMP).
Then, we mark as \textit{troubled} those cells where the candidate DG solution cannot be accepted.  
For these troubled cells we now repeat the time step using, instead of the DG scheme, 
a more robust second order accurate TVD finite volume method, which we assume to produce always an
acceptable solution.

Moreover, in order to maintain the accurate resolution of our original high order DG scheme, 
which would be lost when passing to a FV scheme, 
the FV scheme is applied on a \textit{finer sub-cell grid}, see~\cite{Dumbser2014}.
Indeed, for any troubled cell we define the corresponding sub-cell average of the DG solution at time $t^n$  
\begin{equation}
\label{eq.projection}
\mathbf{v}_{ijk,\alpha}^n(\x,t^n) = \frac{1}{|\omega_{ijk,\alpha}|} \int_{\omega_{ijk,\alpha}^n} \mathbf{u}_{h}^n(\x,t^n) \, d\x 
:=\mathcal{P}(\mathbf{u}_h^n) \qquad \forall \alpha \in [1,N_\omega^d],
\end{equation}
where $|\omega_{ijk,\alpha}|$ denotes the volume of sub-cell $\omega_{ijk,\alpha}$ of element $\Omega_{ijk}$ and $\mathcal{P}(\mathbf{u}_h)$ is the $L_2$ projection operator. 
We then apply a second order TVD FV method in order to evolve $\mathbf{v}_{ijk,\alpha}^n$ and we obtain the FV solution at the sub-cell level at the next time step $\mathbf{v}_{ijk,\alpha}^{n+1}$.
Next, the DG polynomial $\mathbf{u}_{h}^{n+1}$ for each $\Omega_{ijk}$ is recovered from $\mathbf{v}_{ijk,\alpha}^{n+1}$ 
by applying a reconstruction operator $\mathcal{R}$ such that 
\begin{equation}
\int_{\omega_{ijk,\alpha}^n} \mathbf{u}_{h}^{n+1}(\x,t^{n+1}) \, d\x = \int_{\omega_{ijk,\alpha}^n} \mathbf{v}_{ijk,\alpha}^{n+1}(\x,t^n) \, d\x :=\mathcal{R}(\mathbf{v}_{ijk,\alpha}^{n+1}(\x,t^n))  \qquad \forall \alpha \in [1,N_\omega^d],
\label{eq.reconstruction}
\end{equation}
which is \textit{conservative} on the main cell $\Omega_{ijk}$ thanks to the additional linear constraint
\begin{equation}
\int_{\Omega_{ijk}} \mathbf{u}_{h}^{n+1}(\x,t^{n+1}) \, d\x = \int_{\Omega_{ijk}} \mathbf{v}_{h}^{n+1}(\x,t^{n+1}) \, d\x.
\label{eqn.LSQ}
\end{equation}
Finally, note that for the sub-cell FV scheme we have a different CFL stability condition
\be
\label{eq.CFLFV}
\Delta t_{\text{FV}} < \text{CFL}\frac{h_{\text{min}}}{d\, N_\omega}\frac{1}{|\lambda_{\text{max}}|},
\ee
with $h_{\min}$ the minimum cell size referred to $\Omega_{ijk}$.
Condition~\eqref{eq.CFLFV} guides us in choosing the number of employed sub-cells $N_\omega$, and in particular, 
following~\cite{Dumbser2014}, we take $N_\omega = (2N + 1)$ so that $\Delta t_{\text{FV}} = \Delta
t_{\text{DG}}$. This choice allows us to maximize the resolution of the sub-cell FV scheme and to run it at its maximum possible CFL number. 

We conclude this Section with two brief operational remarks on our limiting strategy.
First, note that the reconstruction operator~\eqref{eq.reconstruction}-\eqref{eqn.LSQ} might still lead to an oscillatory solution, 
since it is based on a linear unlimited least squares technique. 
If this is the case, the cell $\Omega_{ijk}$ will be detected as troubled at the next time level $t^{n+2}$, 
therefore the FV sub-cell limiter will be used again in that cell. 
Moreover, in order to overcome the possible issues due to the projection of a non valid reconstructed solution,
the sub-cell averages $\mathbf{v}_{ijk,\alpha}^{n+1}$ are always kept in memory to be reused (instead of recomputed) 
if a cell is detected to be troubled for the second consecutive time step.

Second, in order to keep our scheme \textit{conservative} 
we also need to recompute the DG solution in the non-troubled neighbors (call one of them ${i}$) of a troubled cell (call it ${j}$).
Otherwise at the common space--time lateral surface $\partial \Omega$, 
the flux computed from $\Omega_{{i}}$ would be obtained through the DG scheme, 
while the one coming from the troubled neighbor $\Omega_{{j}}$ would be updated using the sub-cell FV scheme.
Thus, the DG solution in these cells is recomputed in a \textit{mixed way}: 
the volume integral and the surface integrals on good faces are kept, 
while the numerical flux across the troubled faces is always provided by the FV scheme.

\section{Numerical examples}
\label{sec.results}

In all the following numerical examples, the fluid under consideration is described via the ideal gas equation of state with $\gamma_g = 1.4$ and $\pi_g = 0$. \textcolor{black}{To avoid division by zero, the volume fraction $\alpha$ is set in all tests to be in the interval $\alpha \in [\epsilon, 1-\epsilon]$, with $\epsilon$ a
small parameter of the order $10^{-3}$ to $10^{-2}$. In order to allow even smaller values of $\epsilon$
one could use the filtering technique detailed in \cite{Tavelli2019}. } 

\subsection{1D Riemann problems}
\label{sec.rp1d}

The aim of this first series of numerical tests is to verify numerically that the fluid velocity at the 
material interface is indeed equal to the solid velocity when we have a jump in the fluid volume fraction 
from unity to zero. The two-dimensional computational domain is chosen as $\Omega = [-1,+1] \times [-0.1,+0.1]$, 
and is discretized with $100 \times 10$ equidistant ADER-DG elements of polynomial approximation degree $N=3$
and with \textit{a posteriori} sub-cell finite volume limiter. 
The fluid volume fraction $\alpha$ is initially set to $\alpha_L = \epsilon$ for $x \leq 0$ and to $\alpha_R = 1-\epsilon$ for $x > 0$, with $\epsilon = 10^{-3}$. All other quantities are simply initialized with a \textit{constant} value throughout the entire computational domain, hence the entire flow field is  generated  by the jump in the scalar volume fraction function $\alpha$ alone. We consider three scenarios, RP1, RP2 and RP3, with initial data summarized in Table \ref{tab.ic.rp}. 

\begin{table}[!t] 
	 \caption{Initial data for fluid density, fluid velocity, fluid pressure and solid velocity for Riemann problems RP1, RP2 and RP3. Also the final times for each simulation are given.} 
	 \begin{center} 
	 \begin{tabular}{cccccccc}
	 	\hline
	 	    &   $\rho$  & $u$ & $v$ & $p$ & $u_s$ & $v_s$ & $t_{\textnormal{end}}$   \\ 
	 	\hline 
	 	RP1 &   1.0     &  0.0 & 0.0 & 1.0 &  1.0   & 0.0 & 0.4   \\ 
	 	RP2 &   1.0     &  0.0 & 0.0 & 1.0 & -1.0   & 0.0 & 0.4   \\ 
	 	RP3 &   1.0     & -1.0 & 0.0 & 1.0 &  3.0   & 0.0 & 0.2   \\ 
	 	\hline 
	\end{tabular}  	
	 \end{center} 
\label{tab.ic.rp} 
\end{table} 

The first Riemann problem (RP1) represents a solid piston that is moving into a fluid at rest with moderate positive speed, hence causing a right-moving shock wave, while the second Riemann problem (RP2) models a solid piston that is moving away from a fluid at rest with moderate negative speed, hence causing a right-moving rarefaction in the fluid. The last Riemann problem (RP3) describes a piston that hits a left-moving fluid with supersonic velocity and thus generates a rather strong shock, with a shock Mach number of about $M=5$. 

The exact solution of all Riemann problems can be easily found via the exact Riemann solver detailed in the textbook of Toro \cite{toro-book} and making use of Galilean invariance of Newtonian mechanics. 

\begin{figure}[!htbp]
	\begin{center} 
	\begin{tabular}{ccc} 
	\includegraphics[width=0.3\textwidth]{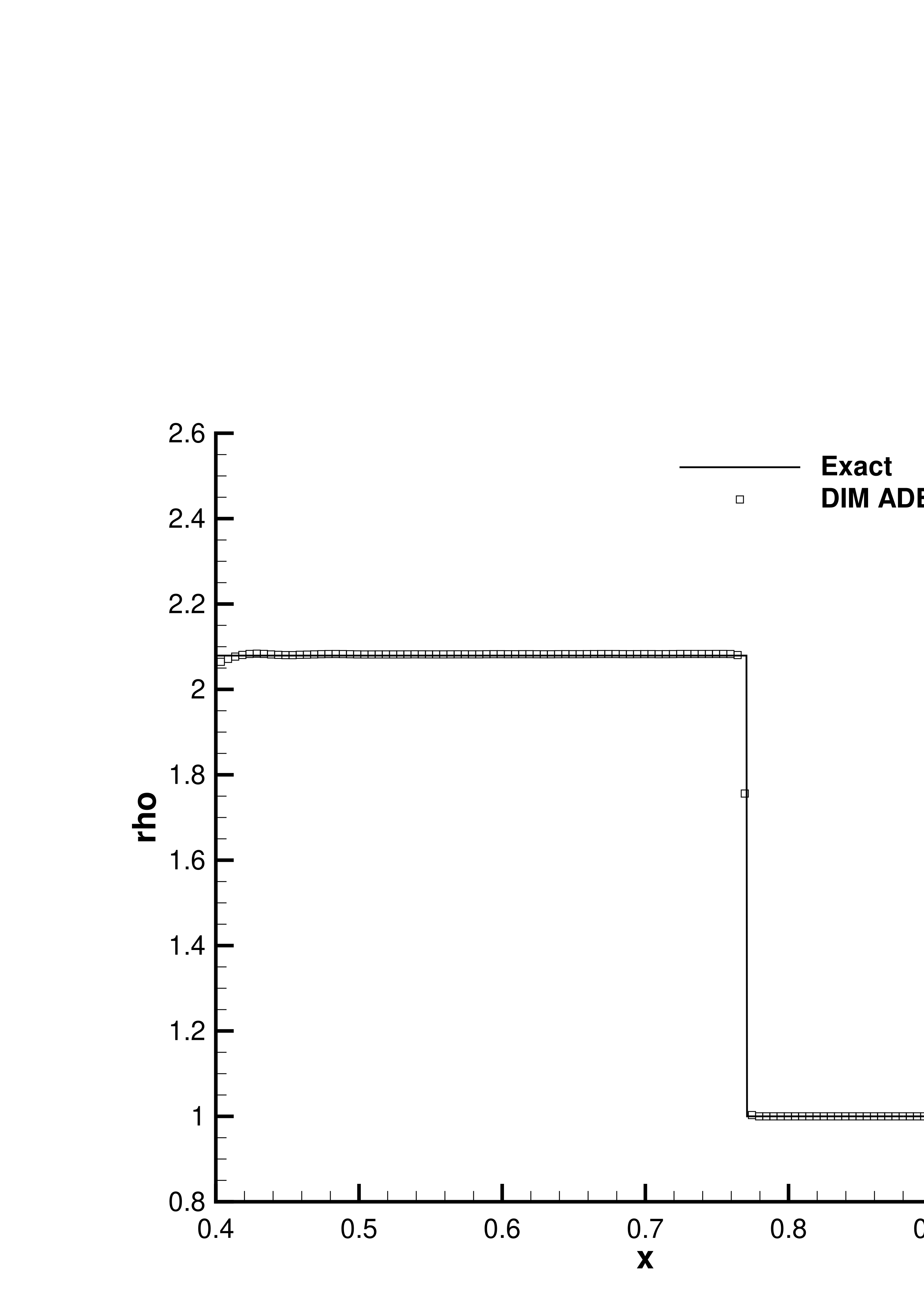}  &
	\includegraphics[width=0.3\textwidth]{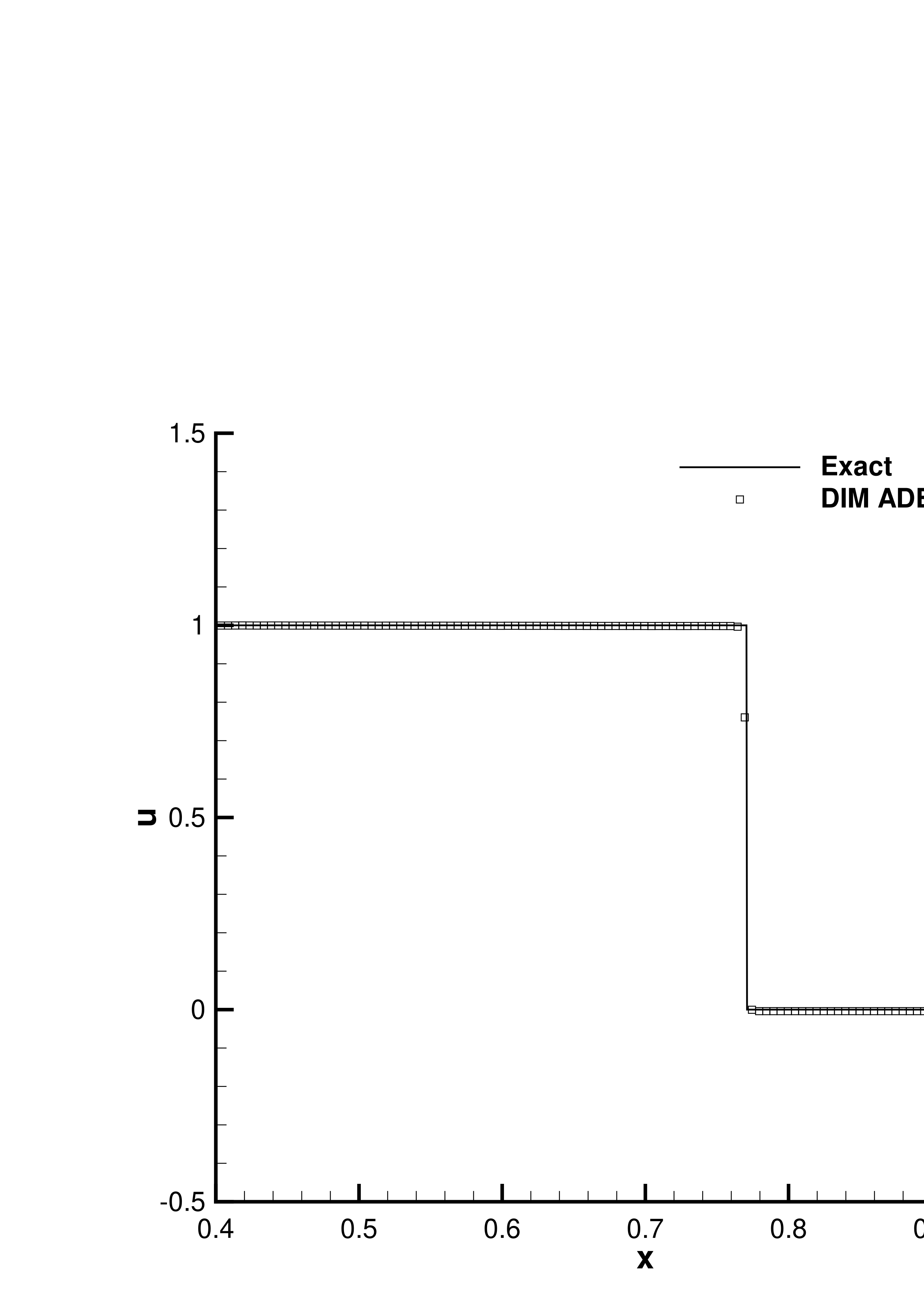}    &
	\includegraphics[width=0.3\textwidth]{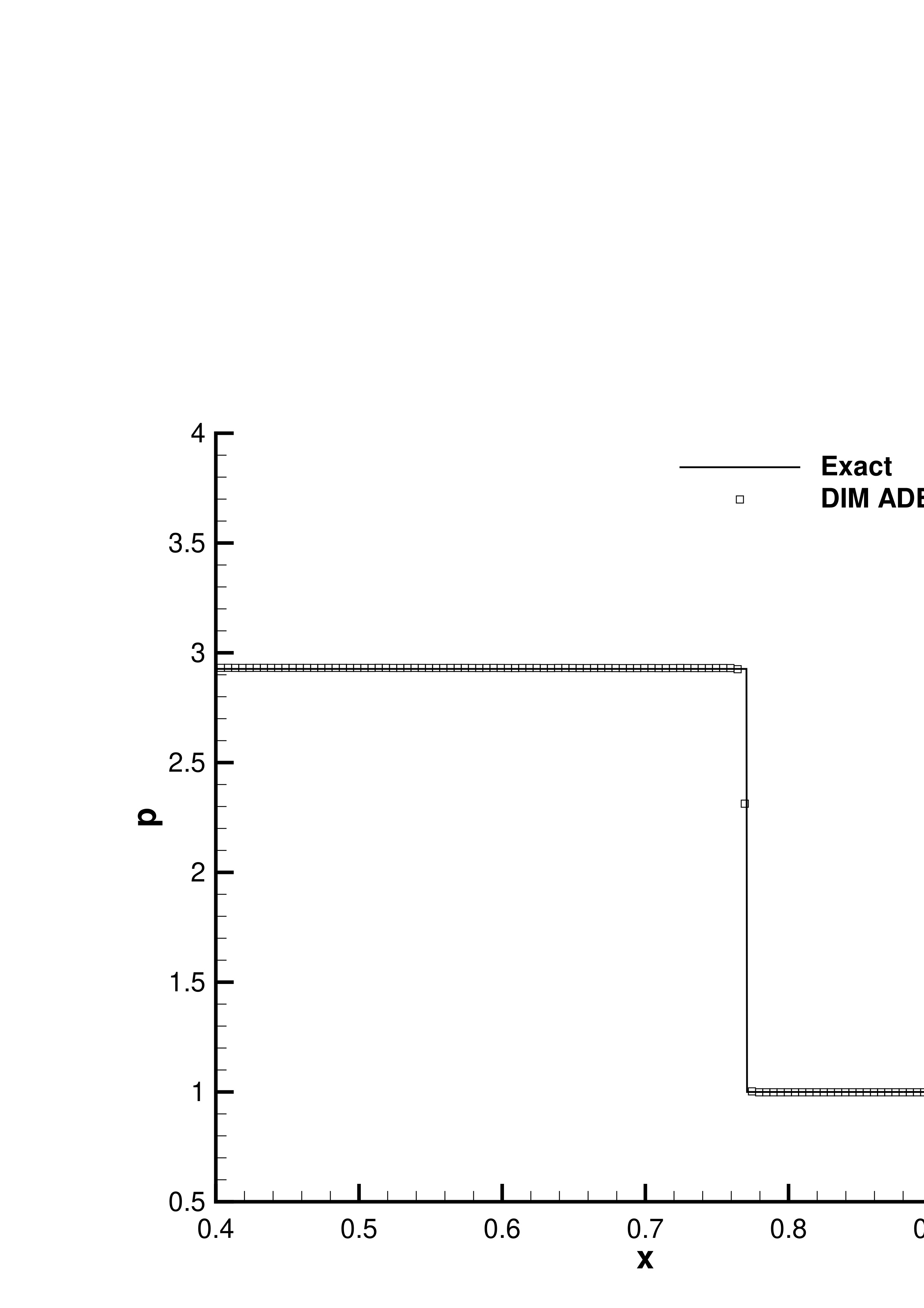}  
	\end{tabular} 
	\end{center} 
	\caption{Riemann problem RP1 (slowly moving piston problem with shock wave) and comparison with the exact solution at time $t=0.4$. Fluid density (left), fluid velocity (center) and fluid pressure (right) in the volume occupied by the fluid. } 
	\label{fig.RP1}
\end{figure}

\begin{figure}[!htbp]
	\begin{center} 
		\begin{tabular}{ccc} 
			\includegraphics[width=0.3\textwidth]{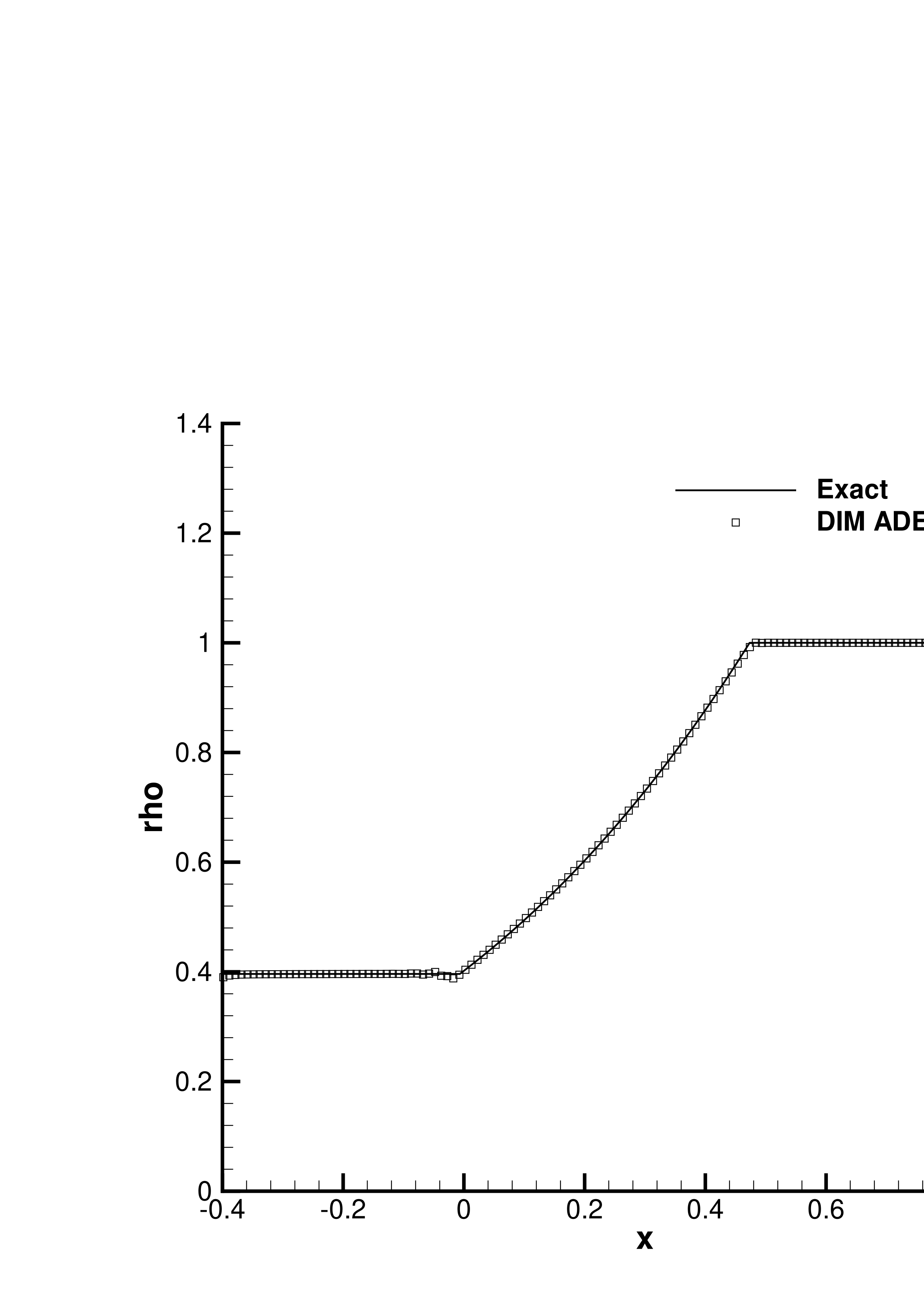}  &
			\includegraphics[width=0.3\textwidth]{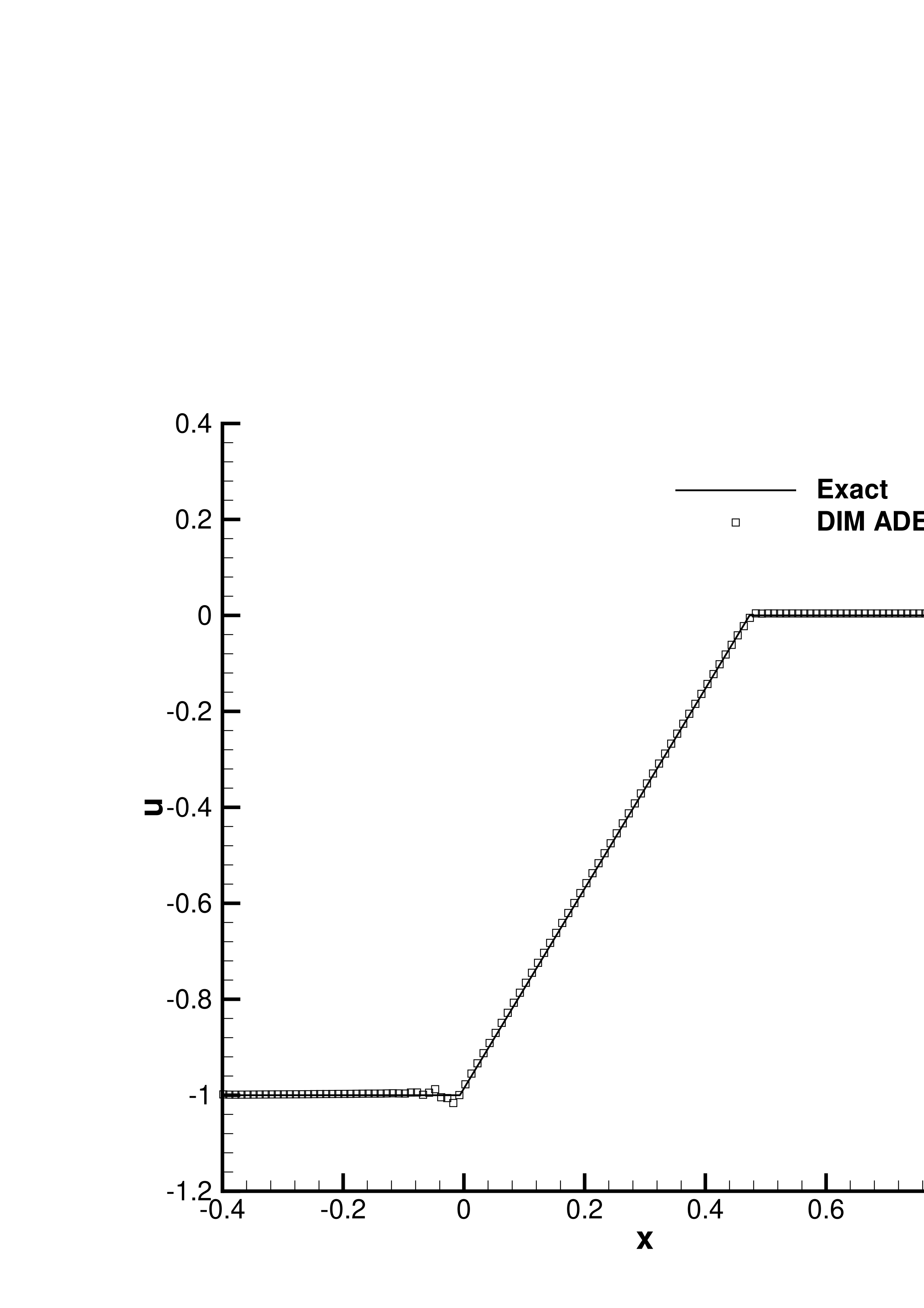}    &
			\includegraphics[width=0.3\textwidth]{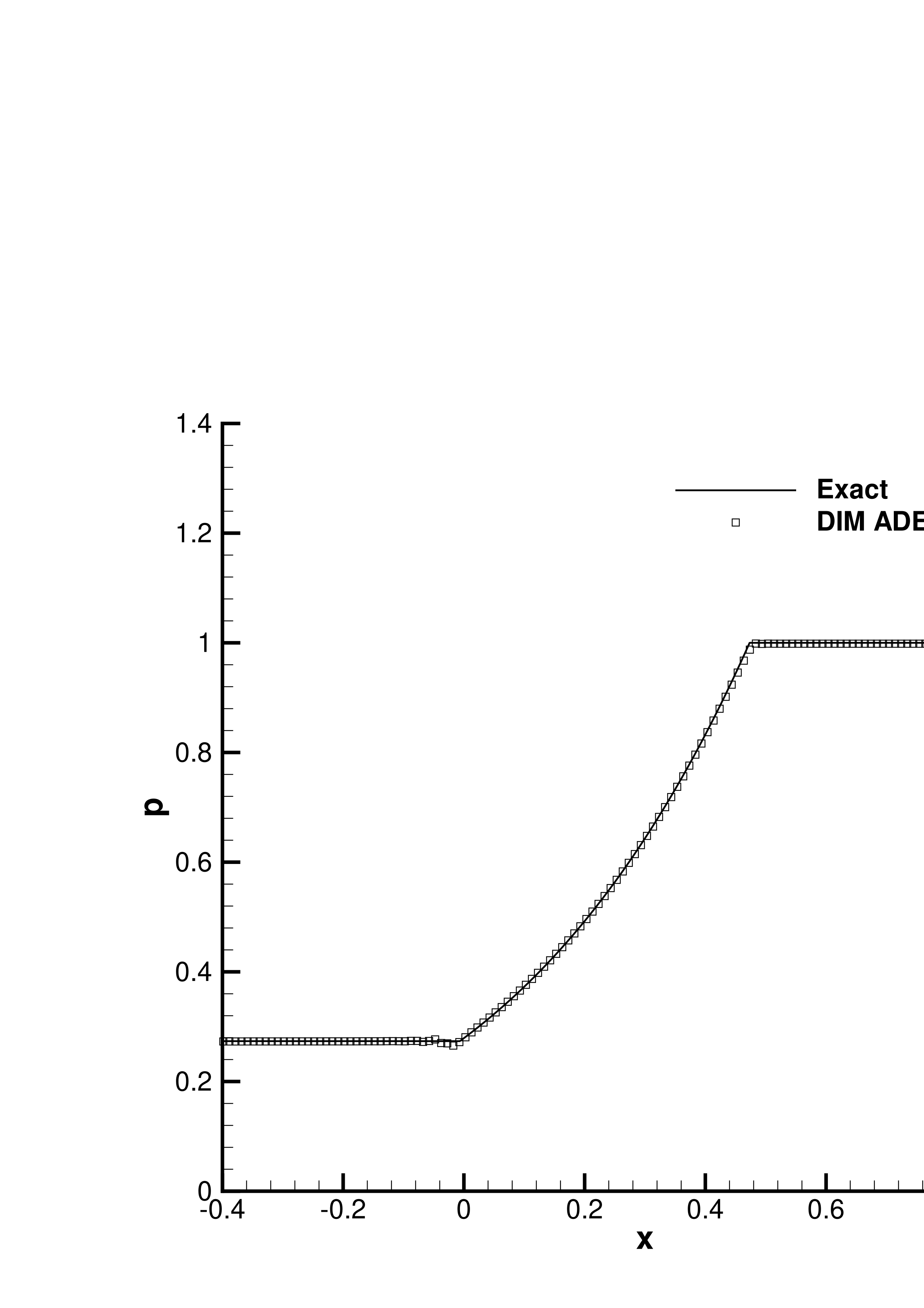}  
		\end{tabular} 
	\end{center} 
	\caption{Riemann problem RP2 (moving piston problem with rarefaction) and comparison with the exact solution at time $t=0.4$. Fluid density (left), fluid velocity (center) and fluid pressure (right) in the volume occupied by the fluid. } 
	\label{fig.RP2}
\end{figure}

\begin{figure}[!htbp]
	\begin{center} 
		\begin{tabular}{ccc} 
			\includegraphics[width=0.3\textwidth]{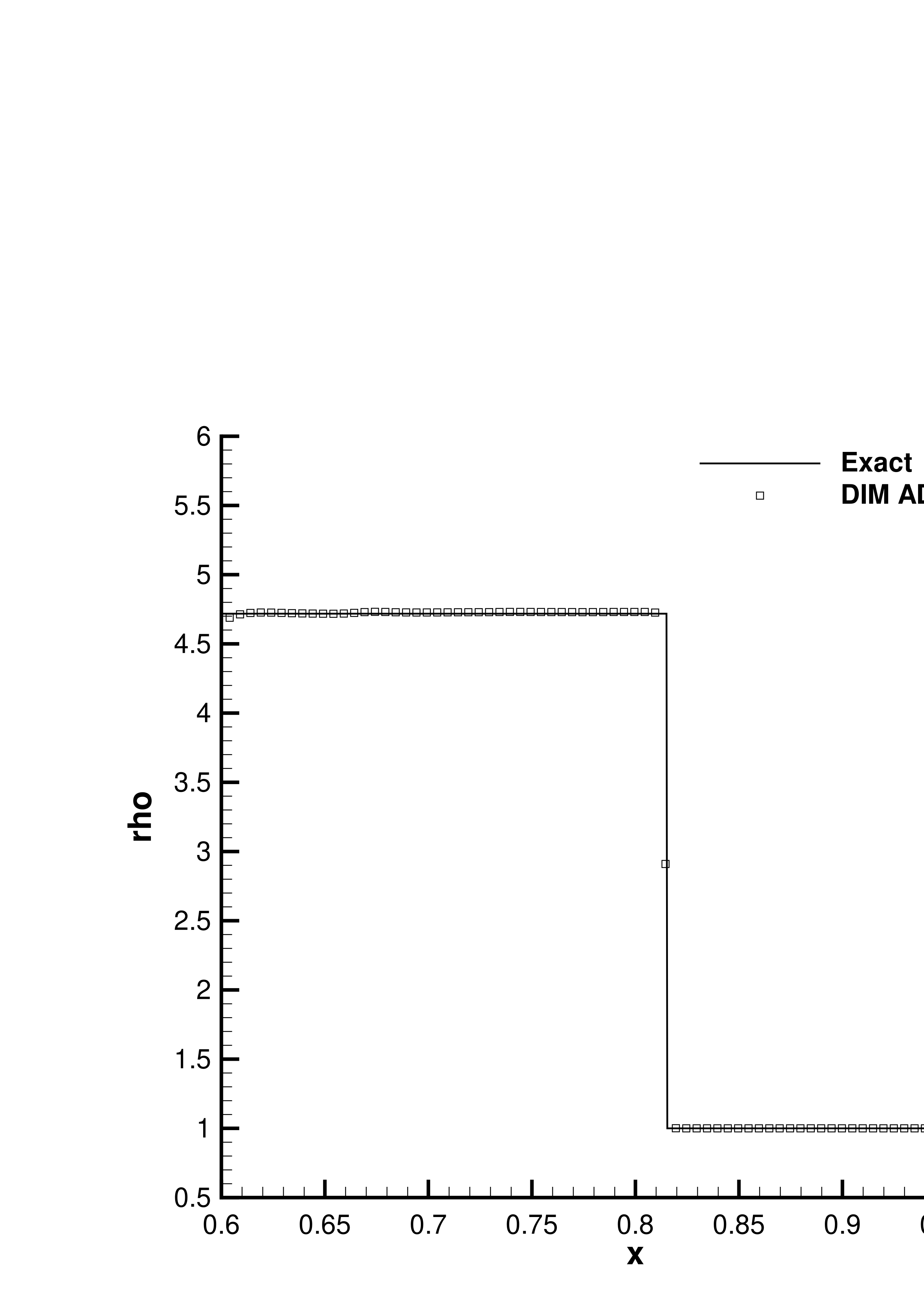}  &
			\includegraphics[width=0.3\textwidth]{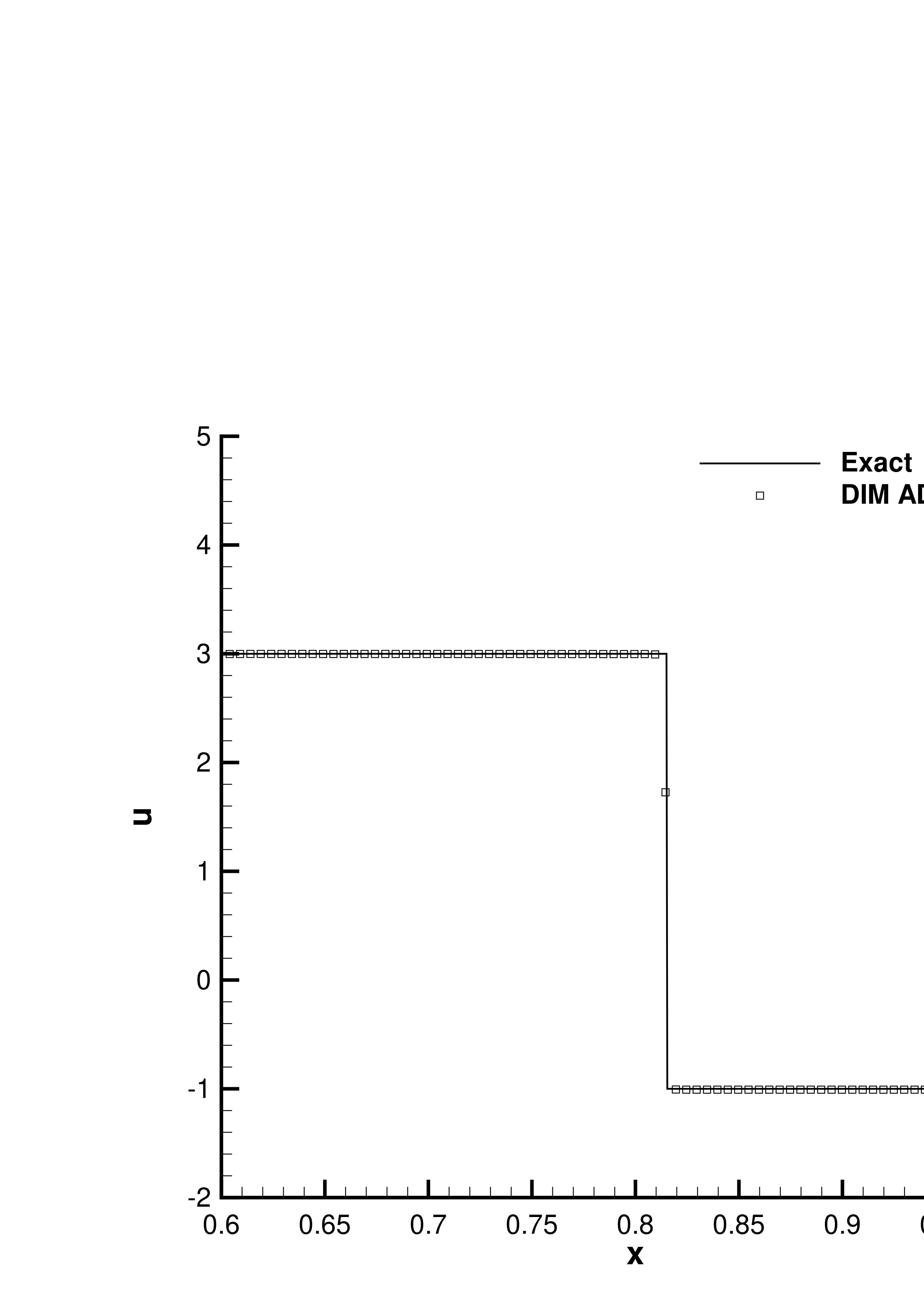}    &
			\includegraphics[width=0.3\textwidth]{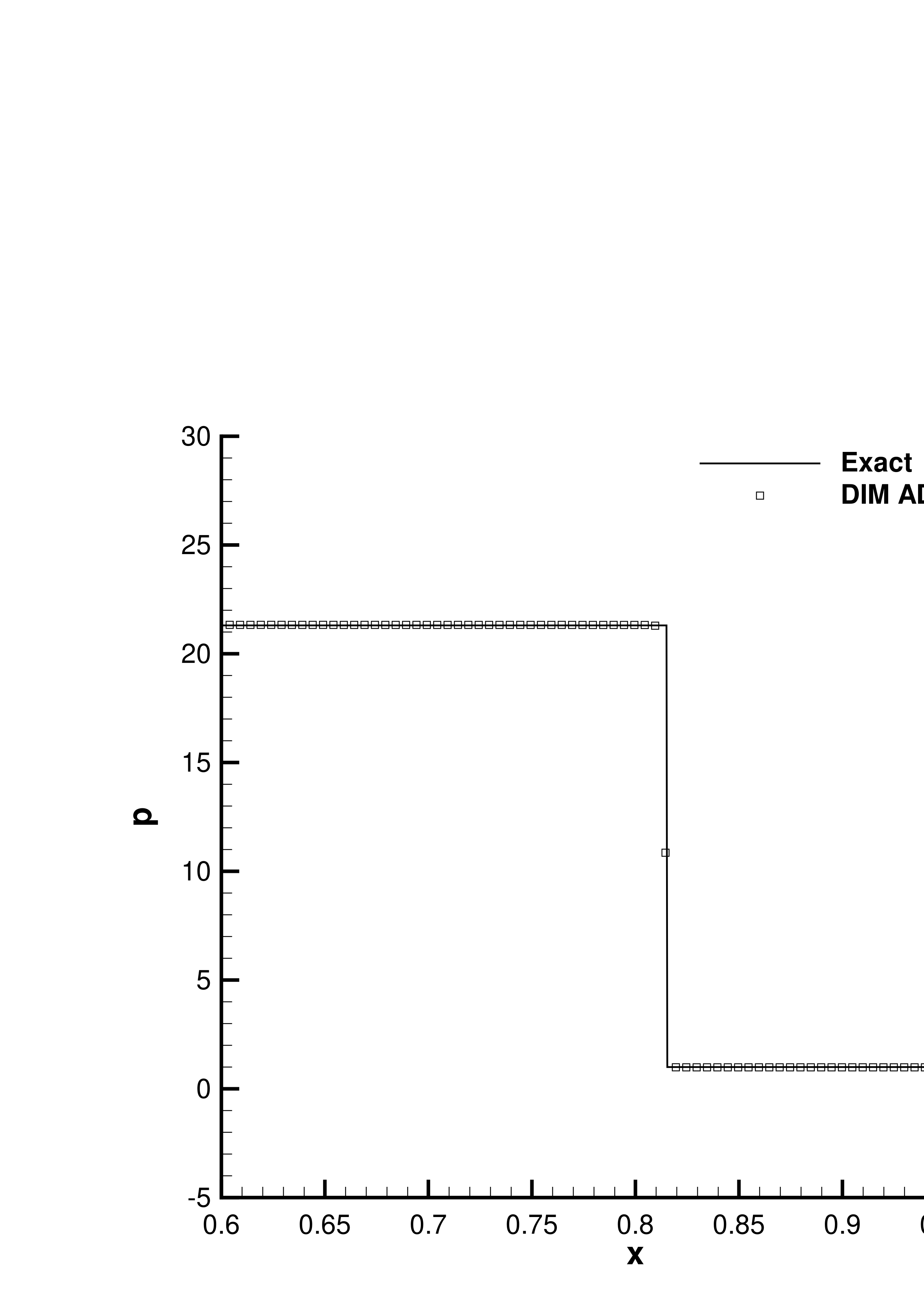}  
		\end{tabular} 
	\end{center} 
	\caption{Riemann problem RP3 (fast moving piston problem with shock wave) and comparison with the exact solution at time $t=0.2$. Fluid density (left), fluid velocity (center) and fluid pressure (right) in the volume occupied by the fluid. } 
	\label{fig.RP3} 
\end{figure}

In Figures \ref{fig.RP1}, \ref{fig.RP2} and \ref{fig.RP3} the numerical results for RP1, RP2 and RP3 are depicted for those parts of the computational domain that are occupied by the fluid at time $t_{\text{end}}$, i.e.
the fluid solid interface is always on the left boundary of each figure.  
In all three cases we can observe an excellent agreement between the numerical solution of the new 
diffuse interface model proposed in this paper and the exact solution. It is also evident that in all simulations the fluid velocity assumes the value of the solid piston, as proven in Section \ref{sec.proof}.

\subsection{Single Mach reflection} 
\label{sec.smr}

Here we solve the singe Mach reflection problem proposed in \cite{toro-book} and for which experimental reference data in form of Schlieren images are available in \cite{toro-book} and \cite{albumfluidmotion}. 
The test consists of a shock wave initially located in $x=0$ and traveling at shock Mach number of $M=1.7$ 
to the right, where it is hitting a wedge of angle $\varphi=25^\circ$. The upstream density and pressure in 
front of the shock ($x>0$) are  $\rho_0=1$ and $p_0 = 1/\gamma_g$, respectively. 
Ahead of the shock, the fluid is at rest ($\mathbf{v}_0=0$). The post-shock values for $x<0$ can then be computed via the standard Rankine-Hugoniot conditions of the compressible Euler equations. 

\begin{figure}[!htbp]
	\begin{center} 
		\begin{tabular}{cc} 
			\includegraphics[width=0.48\textwidth]{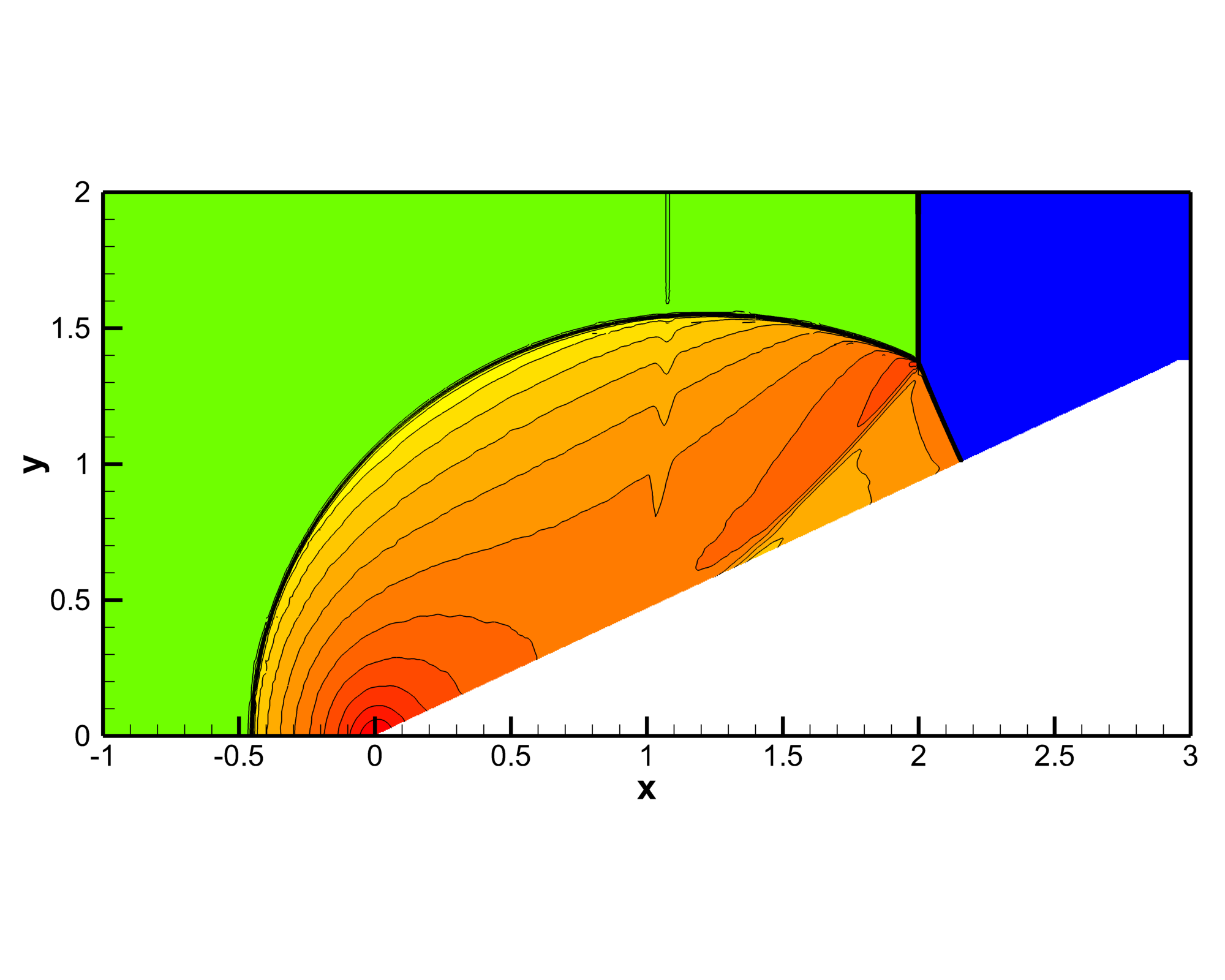}  &
			\includegraphics[width=0.48\textwidth]{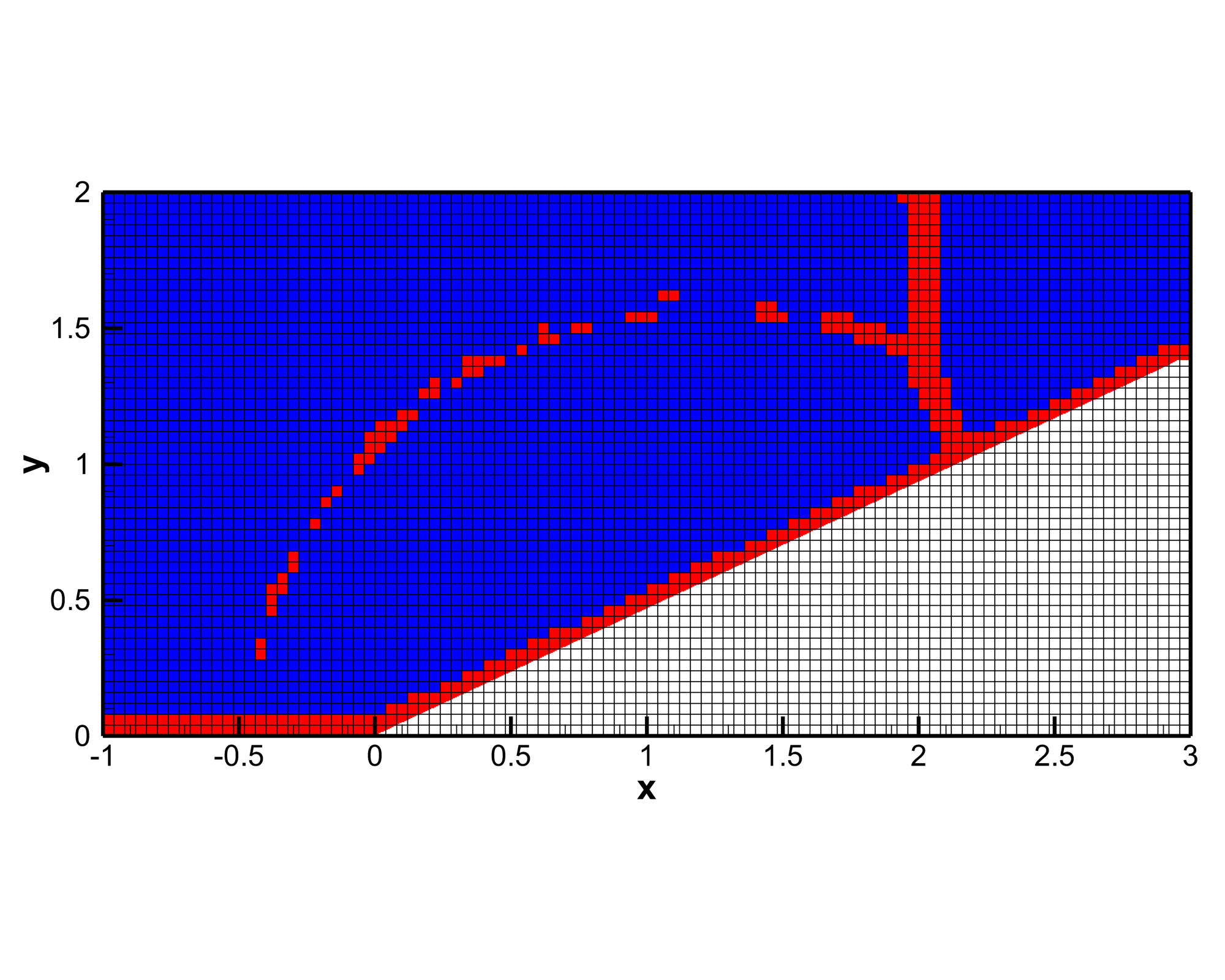} \\ 
			\includegraphics[width=0.48\textwidth]{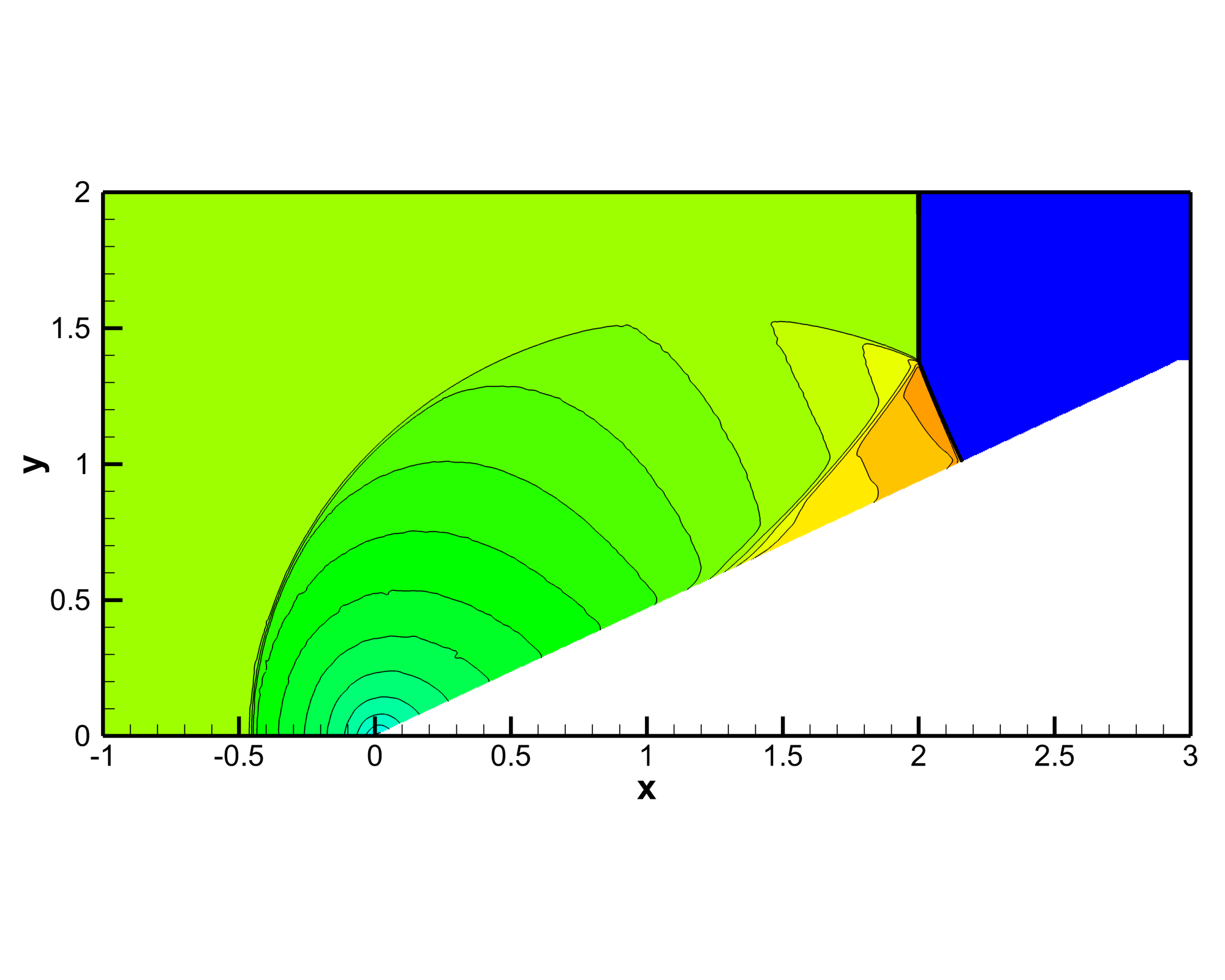}  &
			\includegraphics[width=0.48\textwidth]{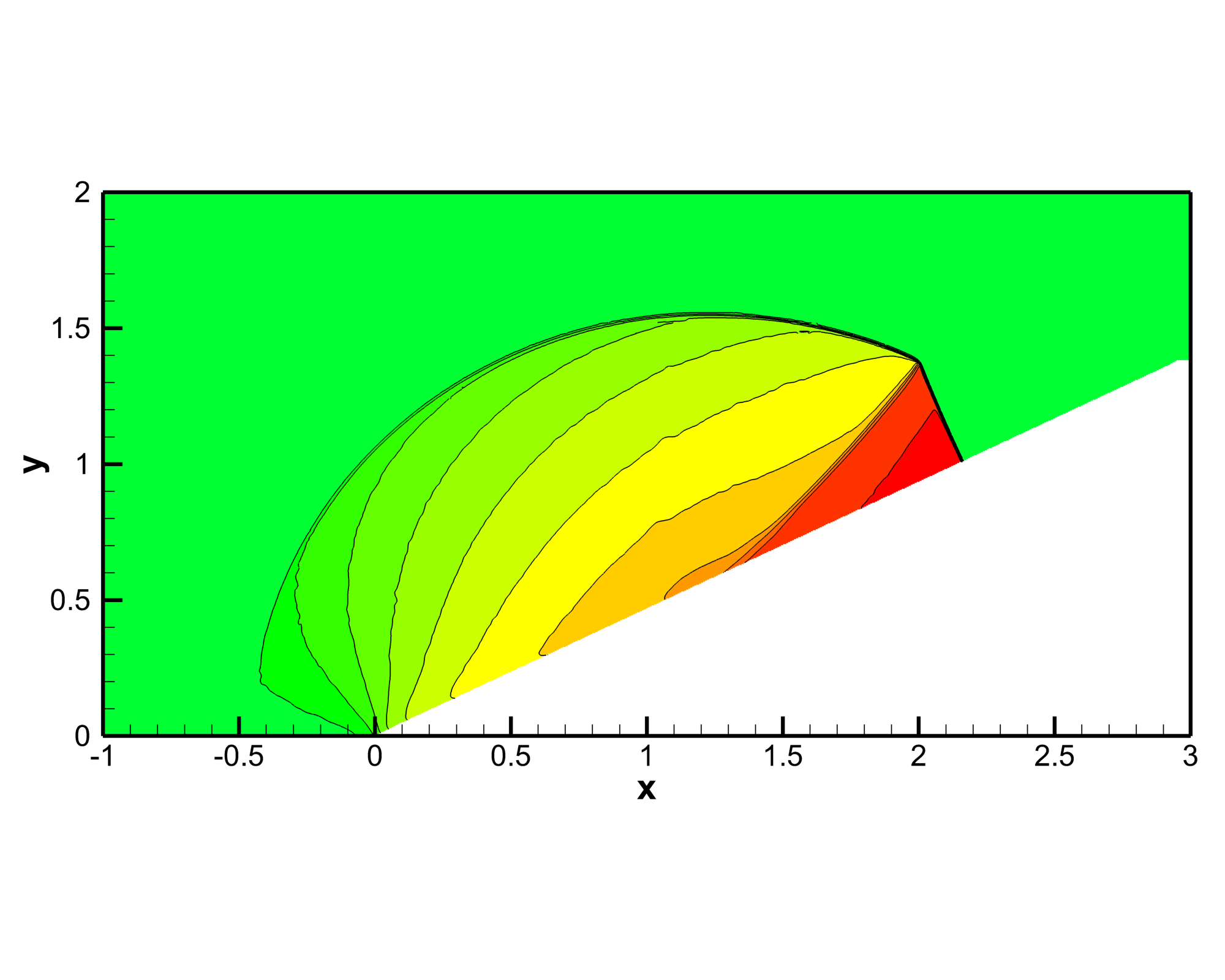}  
		\end{tabular} 
	\end{center} 
	\caption{Single Mach reflection problem at time $t=1.2$ using an ADER-DG scheme with $N=5$ and \textit{a posteriori} sub-cell FV limiter. Density contour colors (top left), limiter map and computational mesh (top right), velocity component $u$ (bottom left) and velocity component $v$ (bottom right).  } 
	\label{fig.smr} 
\end{figure}

The computational domain is $\Omega = [0,3] \times [0,2]$ and is discretized using $100 \times 50$ ADER-DG  elements of degree $N=5$. 
The initial volume fraction function is chosen as $\alpha = \epsilon$ for $y < \tan( \varphi ) \, x$ and
$\alpha = 1-\epsilon$ elsewhere. The velocity of the solid is set to $\mathbf{v}_s=0$. The simulation 
is run until $t=1.2$ and the obtained computational results are summarized in Figure \ref{fig.smr}. 
The obtained flow field agrees very well with the reference solution shown in \cite{toro-book}. The shock
is in the right location and is well resolved with our high order scheme, although only a very coarse mesh
has been used. The limiter is activated along the fluid-solid interface, and along the impinging and reflected shock waves. We emphasize again that all that is needed in order to represent the geometry of the rigid 
solid body is to set the volume fraction function $\alpha$ to a small positive value inside the solid, 
and to almost unity outside.

\subsection{Flow of a shock over a wedge}  
\label{sec.wedge}

Here we solve a similar problem as before by considering the interaction 
of a mild shock wave with a two dimensional wedge. For 
this test, experimental reference data are again available under the form 
of Schlieren photographs, see \cite{albumfluidmotion,schardin,DumbserKaeser07}. 

The computational domain is $\Omega = [-2,6]\times[-3,3]$ and is discretized
with $200 \times 150$ ADER-DG elements with polynomial approximation degree $N=5$. 
The scheme is supplemented with \textit{a posteriori} sub-cell finite volume limiter. 
The tip of a wedge with length $L=1$ and height $H=1$ is placed at 
$x=0$. Inside the wedge we set the initial volume fraction function to $\alpha=\epsilon$, 
and outside we use $\alpha=1-\epsilon$, with $\epsilon = 10^{-2}$. 
In the gas, the exact solution of a shock wave of shock Mach number $M_s=1.3$ is 
setup via the Rankine-Hugoniot conditions of the compressible Euler equations. 
In front of the shock, the fluid is at rest and has a density of $\rho=1.4$ and a
pressure of $p=1$. The solid velocity is set to $\mathbf{v}_s=0$ everywhere. 

\begin{figure}
	\begin{center} 
	\begin{tabular}{cc}
		\includegraphics[width=0.48\textwidth]{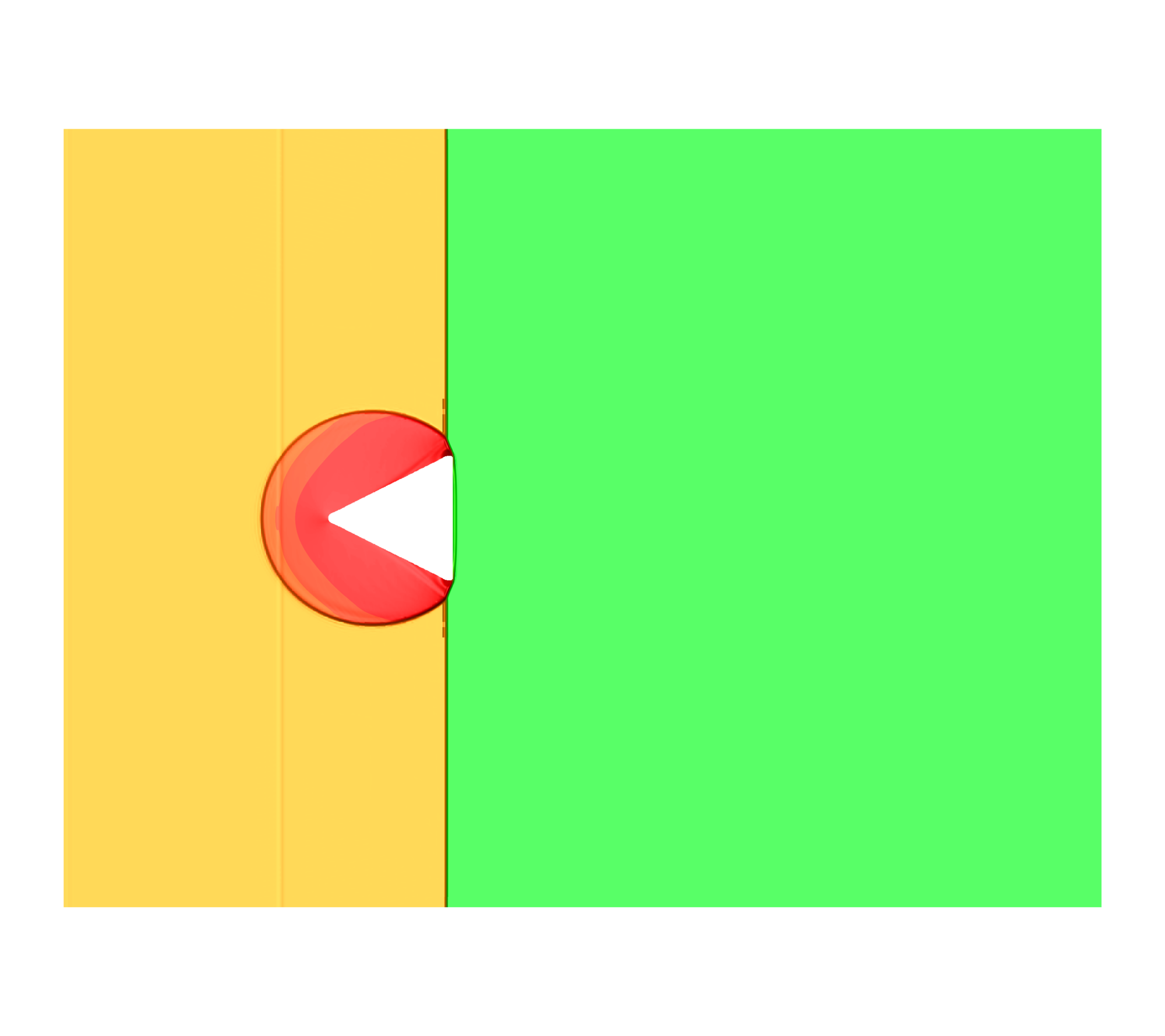}  & 
		\includegraphics[width=0.48\textwidth]{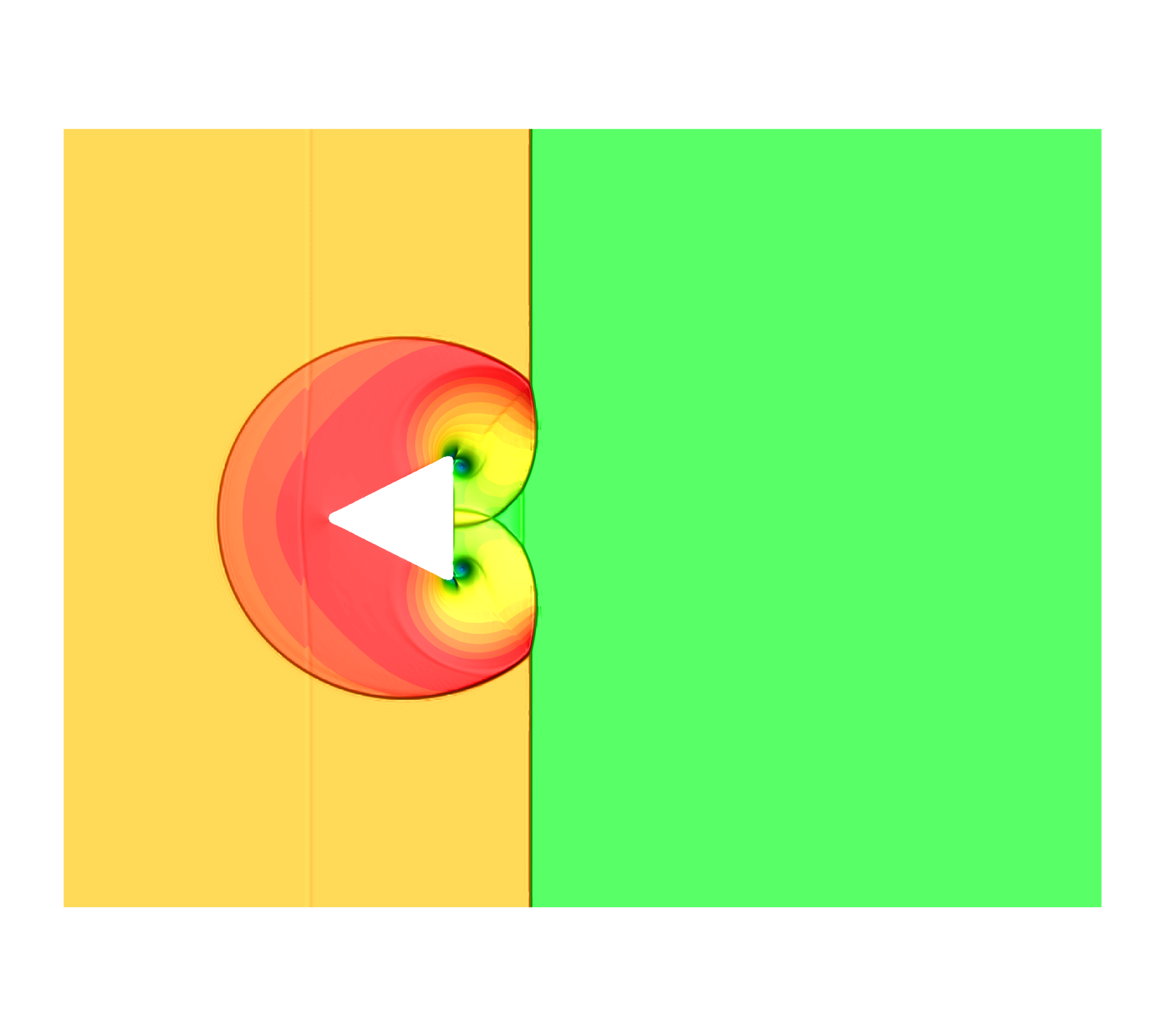}  \\ 
		\includegraphics[width=0.48\textwidth]{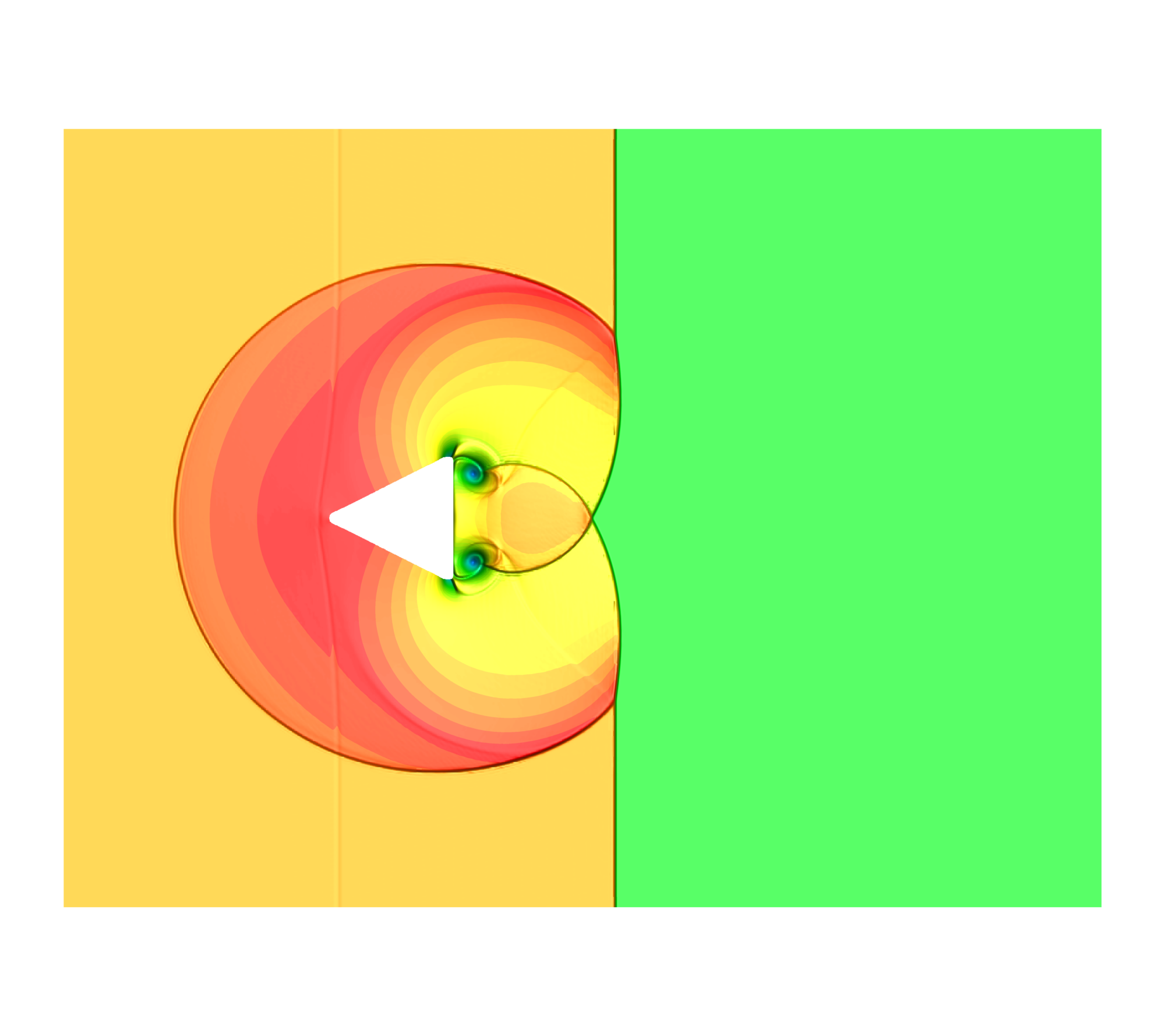}  & 
        \includegraphics[width=0.48\textwidth]{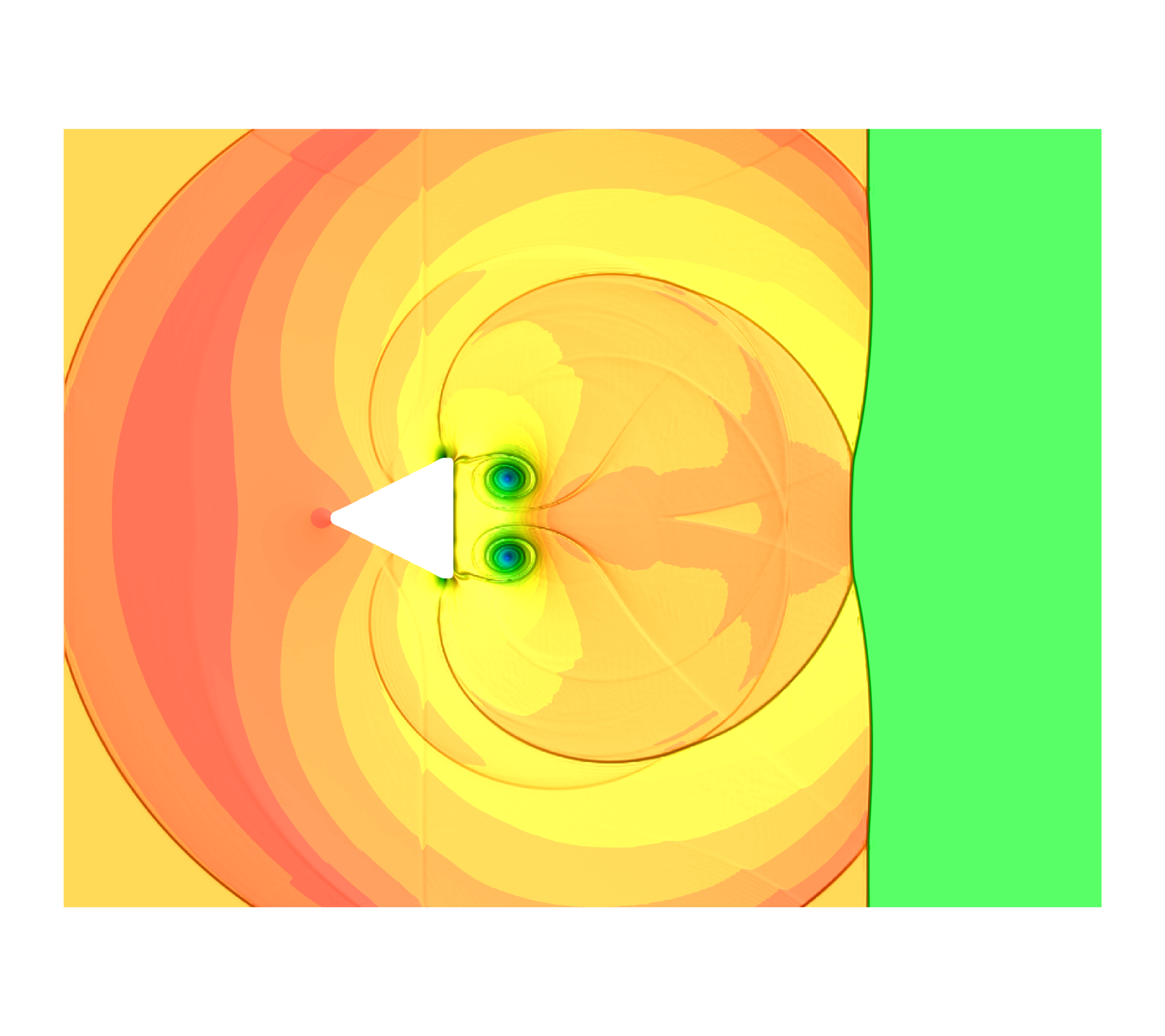}    
	\end{tabular}
	\end{center}
	\caption{Density contours obtained with the diffuse interface model for the shock-wedge 
			interaction problem using a sixth order ADER-DG scheme with \textit{a posteriori} 
			sub-cell finite volume limiter. 
			Output times from top left to bottom right: $t=1.5$, $t=2.0$, $t=2.5$, $t=4.0$.} 
	\label{fig.shockwedge}
\end{figure}

The computational results are depicted in Figure \ref{fig.shockwedge}. Comparing our results 
qualitatively with the Schlieren images produced by Schardin \cite{schardin}, 
we note an excellent agreement. 
All the reflected and refracted waves as well as the vortices shed from the 
wedge are resolved correctly. Our set of images corresponds to pictures 
2, 5, 8 and 17 shown in \cite{schardin}. For a numerical reference solution obtained
with high order ADER-FV schemes on a very fine boundary-fitted unstructured triangular 
mesh, see \cite{DumbserKaeser07}. 

Again, this rather complex flow field is simply produced
via a spatially variable volume fraction function $\alpha$, which defines the geometry
of the obstacle. In principle, with our new approach obstacles of arbitrary shape could
be described.

\begin{figure}
	\begin{center} 
		\begin{tabular}{cc}
			\includegraphics[width=0.48\textwidth]{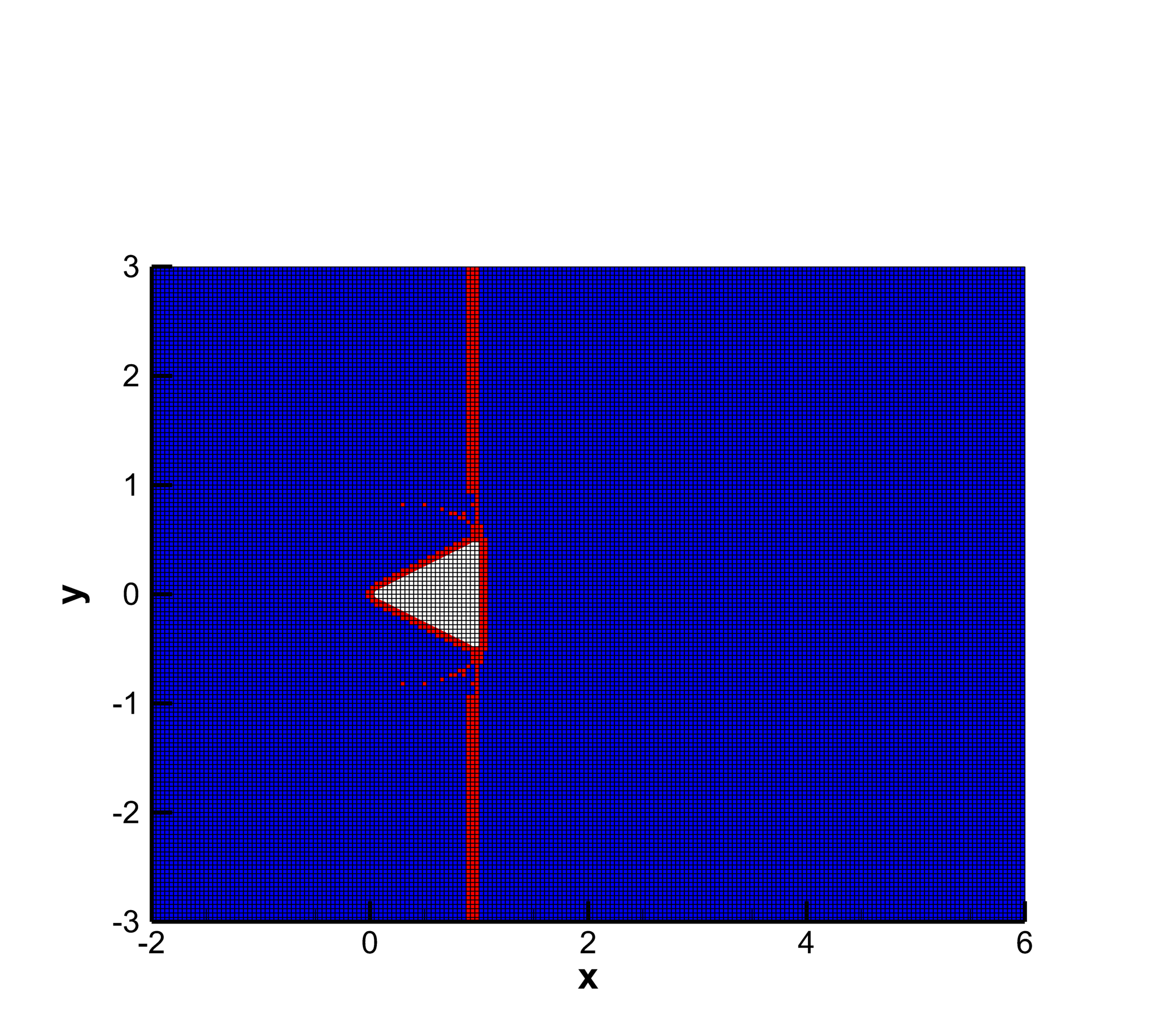}  & 
			\includegraphics[width=0.48\textwidth]{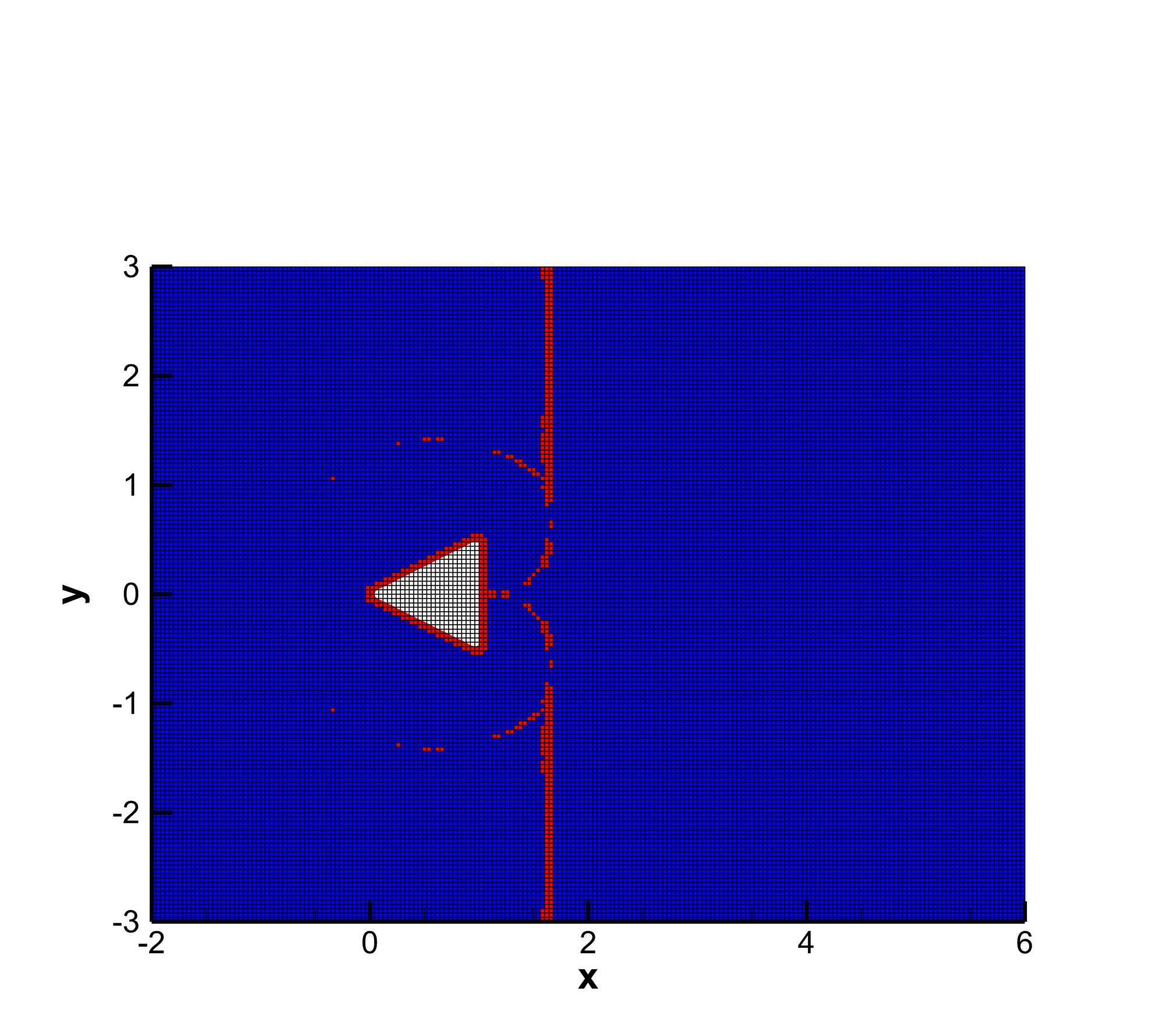}  \\ 
			\includegraphics[width=0.48\textwidth]{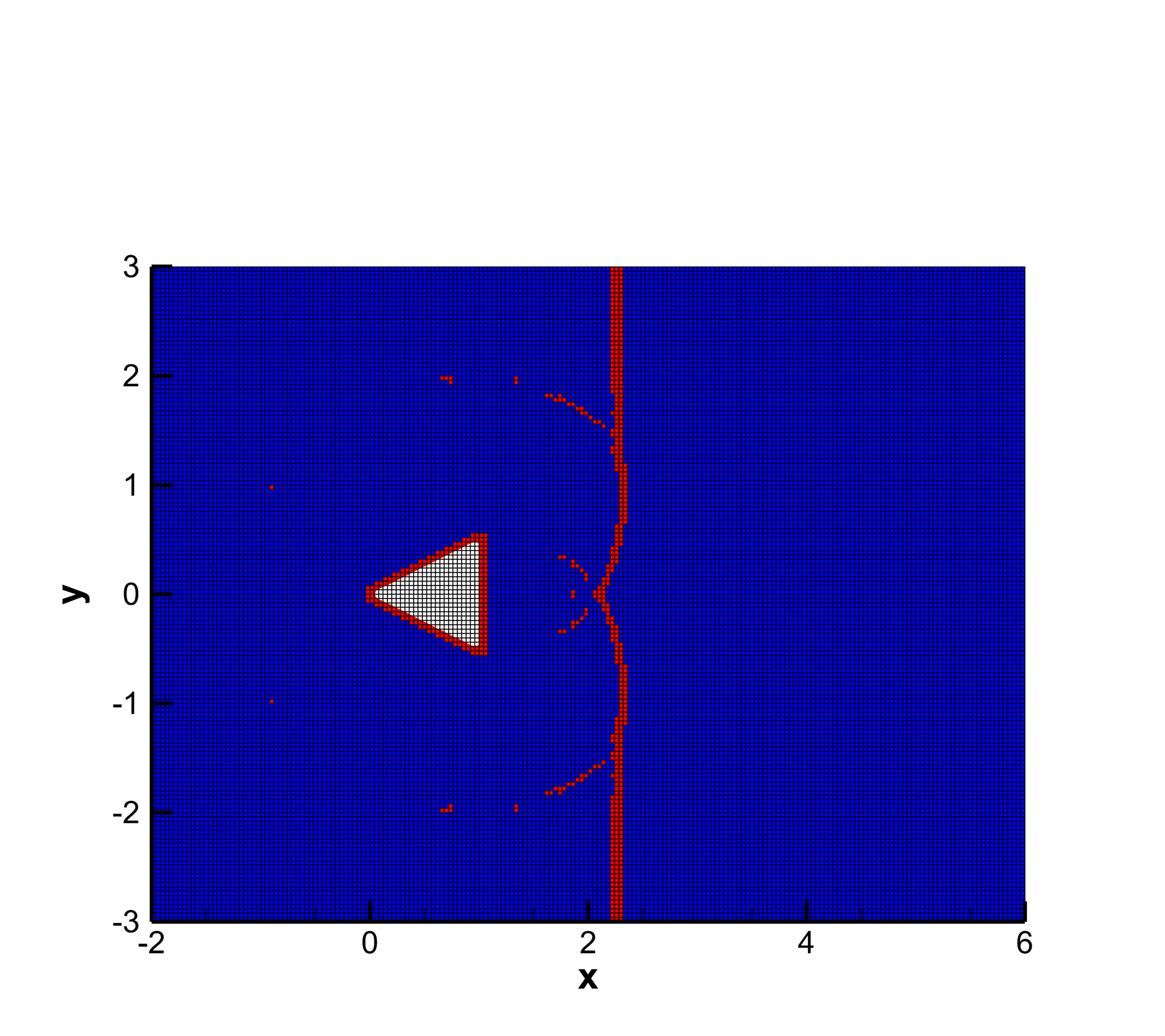}  & 
			\includegraphics[width=0.48\textwidth]{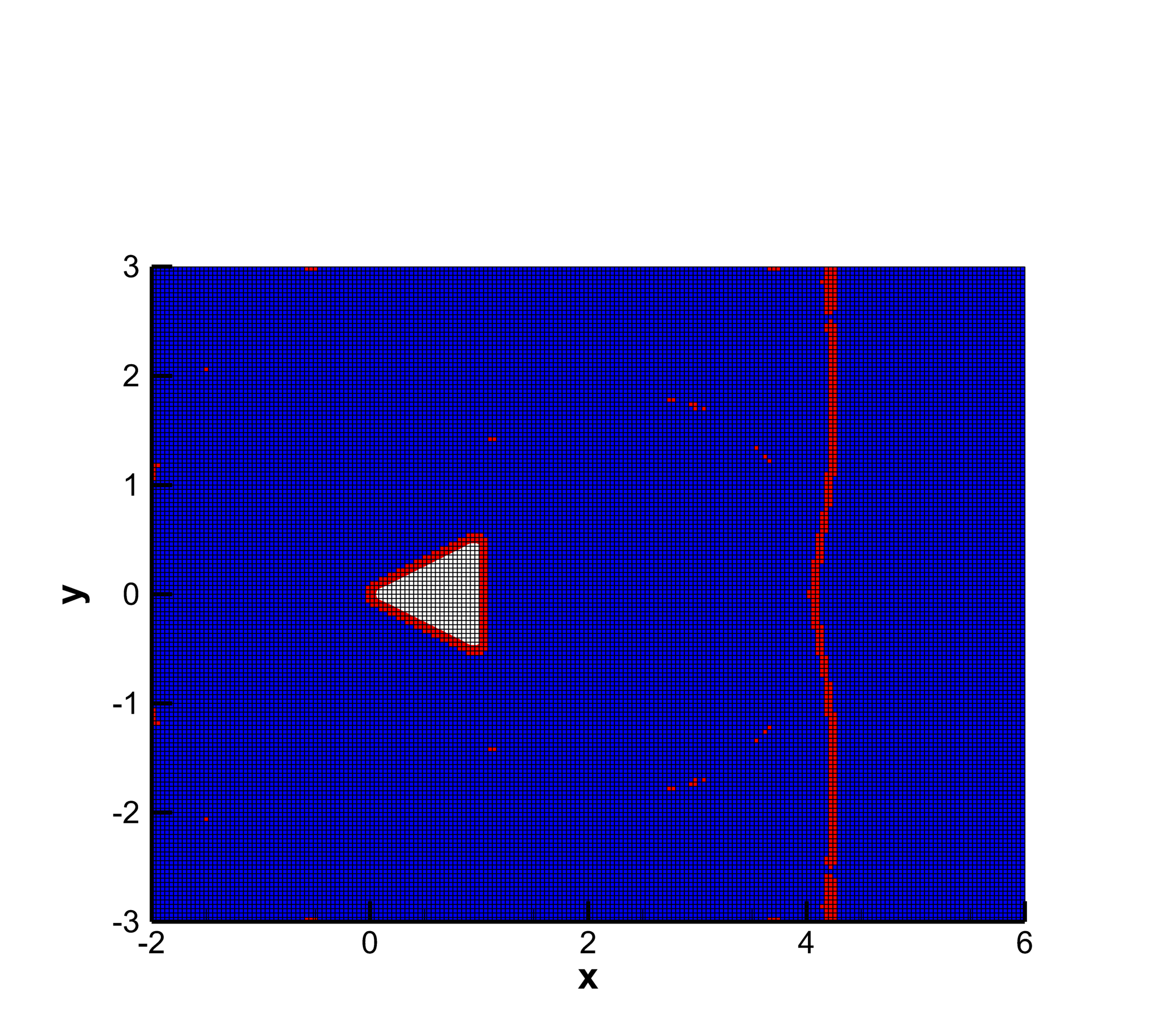}    
		\end{tabular}
	\end{center}
	\caption{Limiter map for the shock-wedge interaction problem using a sixth order 
		ADER-DG scheme with \textit{a posteriori} sub-cell finite volume limiter.  
		Blue cells are unlimited, while red cells have been flagged as troubled. 
		Output times from top left to bottom right: $t=1.5$, $t=2.0$, $t=2.5$, $t=4.0$.} 
	\label{fig.shockwedge.lim}
\end{figure}

The temporal evolution of the limiter is shown in Figure \ref{fig.shockwedge.lim}. The limiter
is activated on the fluid-solid interface and on the main shock wave, as expected.

\subsection{Mach 3 flow over a blunt body} 
\label{sec.blunt}

\begin{figure}[!htbp]
	\begin{center} 
		\begin{tabular}{ccc} 
			\includegraphics[width=0.33\textwidth]{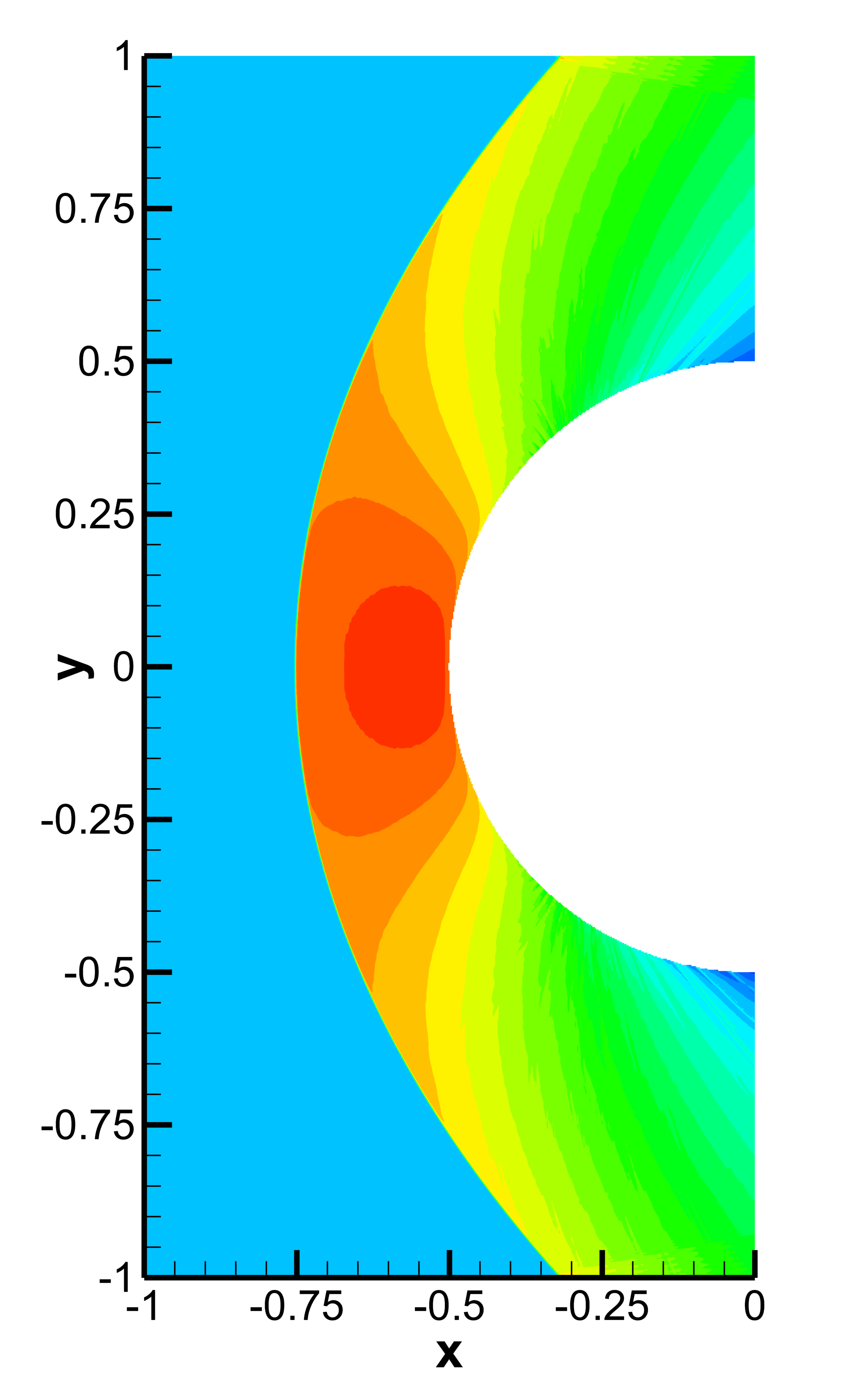}            &
			\includegraphics[width=0.33\textwidth]{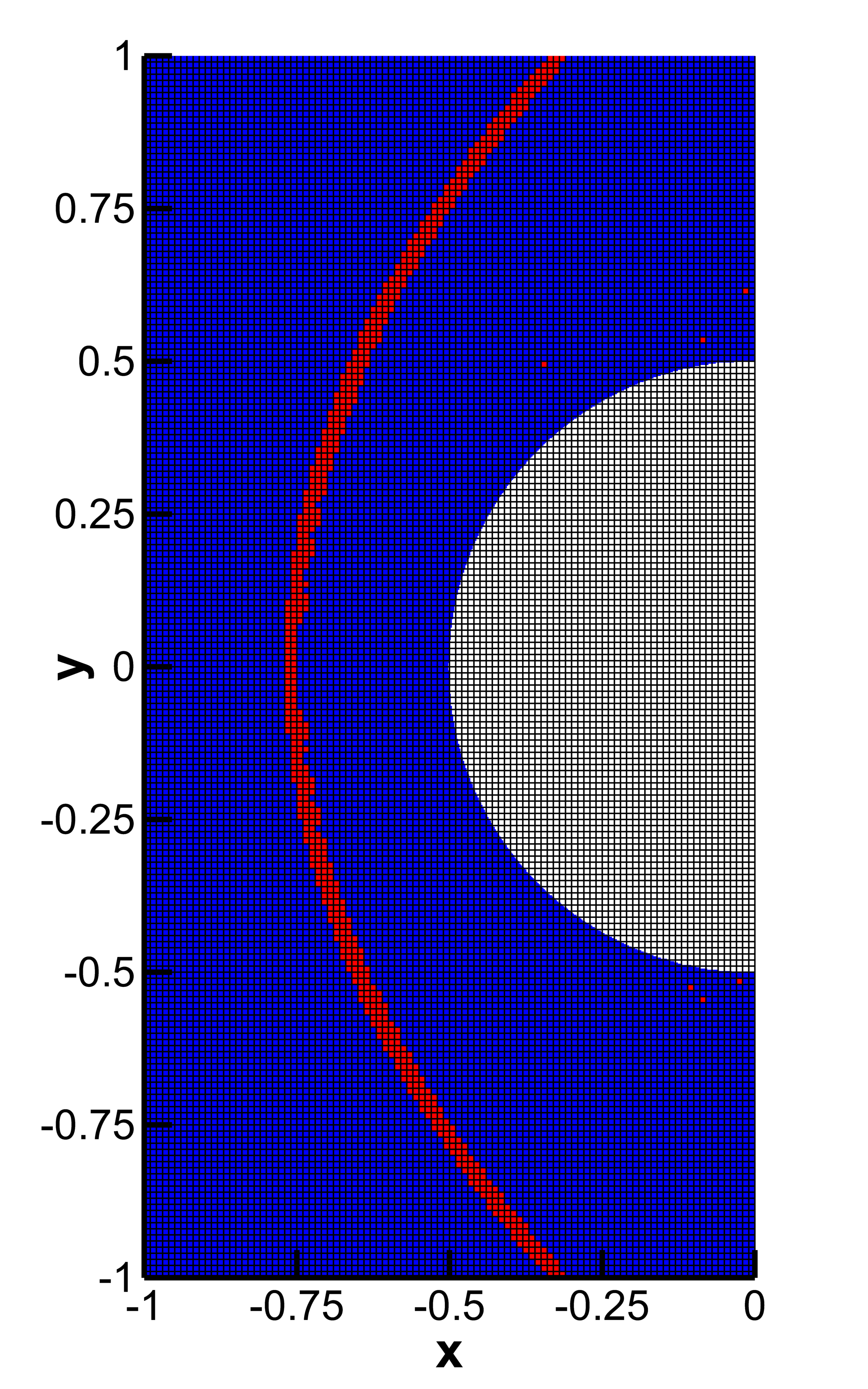}    & 
			\includegraphics[width=0.33\textwidth]{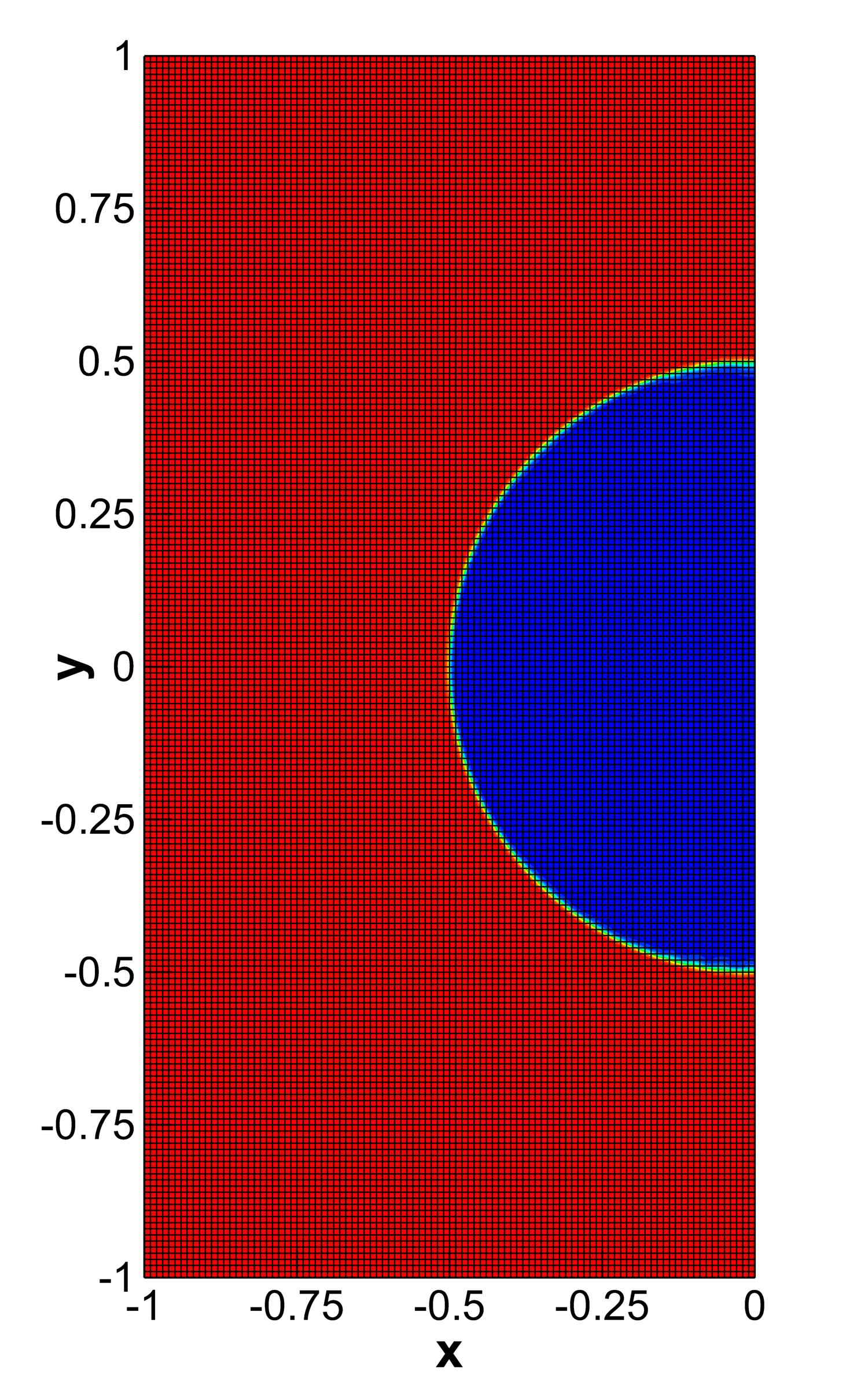}      
		\end{tabular} 
	\end{center} 
	\caption{Supersonic flow over a blunt body with Mach number $M=3$ at time $t=1$ using an ADER-DG scheme with $N=5$ and \textit{a posteriori} sub-cell finite volume limiter. Density contour colors (left) and limiter map with uniform Cartesian grid (center). Only cells with $\alpha > 0.5$ are shown in the left and center panel. Equidistant alpha contour colors (right). The interface is typically resolved within one or two grid cells. } 
	\label{fig.blunt} 
\end{figure}

The next example concerns a supersonic flow with Mach number $M=3$ over a circular cylinder. The two-dimensional computational domain is $\Omega = [-1, 0] \times [-1, +1]$, which is discretized using $100 \times 200$ ADER-DG elements of polynomial approximation degree $N=5$ and \textit{a posteriori} sub-cell finite volume limiter. 
The initial condition for the volume fraction function is simply chosen as follows: inside a circular cylinder of radius $R=0.5$ we set $\alpha=\epsilon$ for $\left \| \mathbf{x} \right\| \leq 0.5$, while outside the cylinder, we set $\alpha = 1 - \epsilon$ in the rest of the computational domain. Here, we use $\epsilon = 10^{-2}$. All other flow quantities are then simply initialized with the following \textit{constant} values: 
$\rho = 1.4$, $u=3$, $v=0$, $p=1$, $u_s=0$ and $v_s=0$. The computational results obtained at
time $t=1$ are depicted in Figure \ref{fig.blunt}. The typical detached bow shock ahead of the blunt body is formed. In the right panel of the same Figure, we also show the limiter map, where red elements are flagged as troubled, while blue cells indicate unlimited elements. One can not only easily see that the \textit{a posteriori} MOOD algorithm is able to detect the shock wave properly, but also that our very high order
scheme is able to resolve the geometry of the blunt body and the shock wave very well, despite using 
a very coarse mesh. \textcolor{black}{Thanks to the subcell FV limiter, the interface can typically be resolved within one to two elements of the DG scheme, see right panel in Figure~\ref{fig.blunt}. In the limiter map one can observe some false positive activations of the limiter in the smooth subsonic area behind the bow shock. However, these local effects do not reduce the overall quality of the simulation. Since the entire scheme is nonlinear and no special care about symmetry preservation was made, such false positive activations can be non-symmetric and can be triggered either via spurious numerical oscillations, accumulated roundoff errors or by passing acoustic waves. The main flow field, however, remains essentially symmetric w.r.t. the $x$ axis. } 

We stress again that also in this test case the entire flow field and the bow shock are merely generated by the 
spatially varying volume fraction function, which is the only information that is needed in order to 
represent the geometry of the solid body. 

\subsection{Three cylinders rotating in a compressible gas at supersonic speed} 
\label{sec.multicyl}

\begin{figure}
	\begin{center} 
		\begin{tabular}{cc}
			\includegraphics[width=0.48\textwidth]{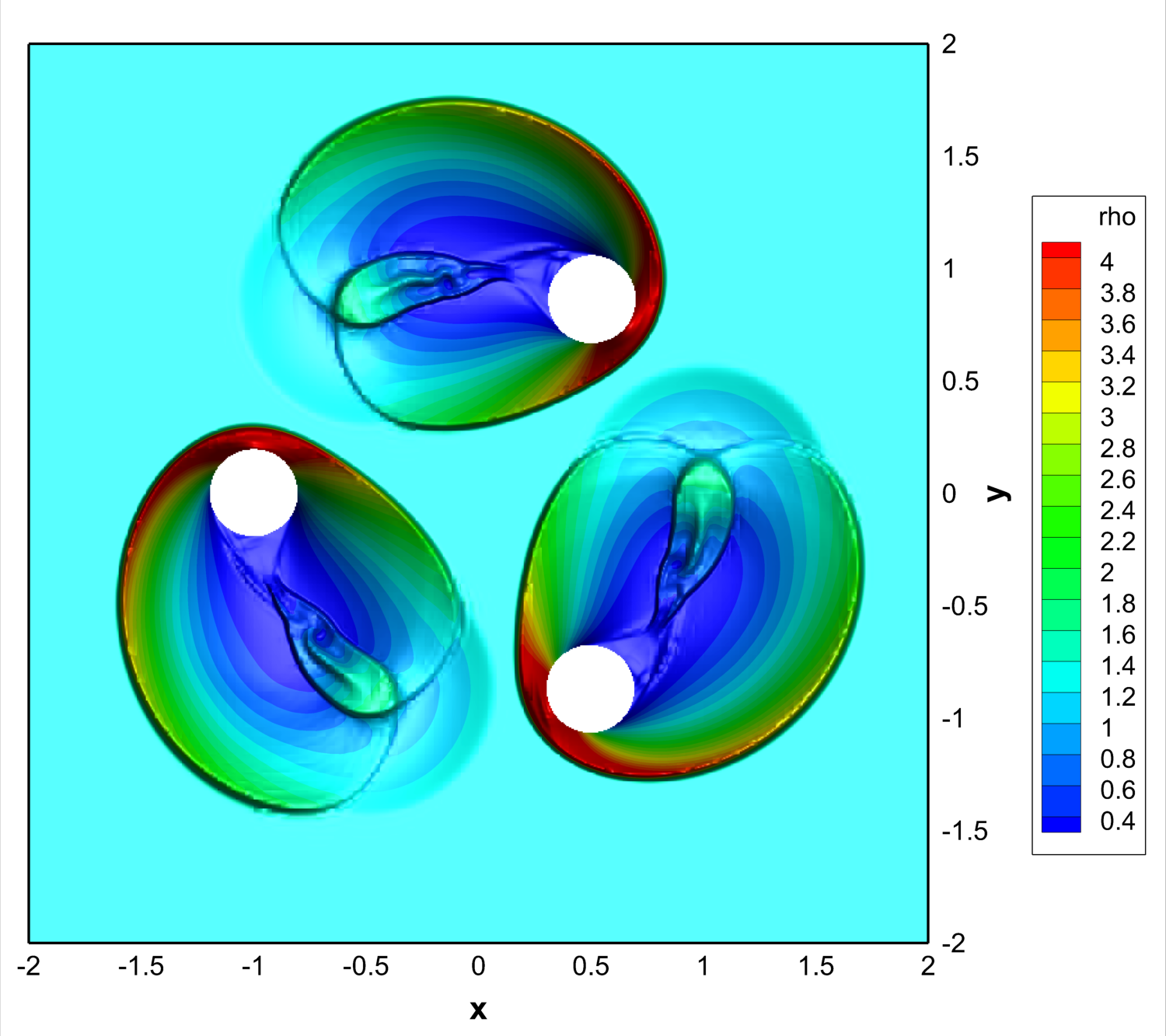}  & 
			\includegraphics[width=0.48\textwidth]{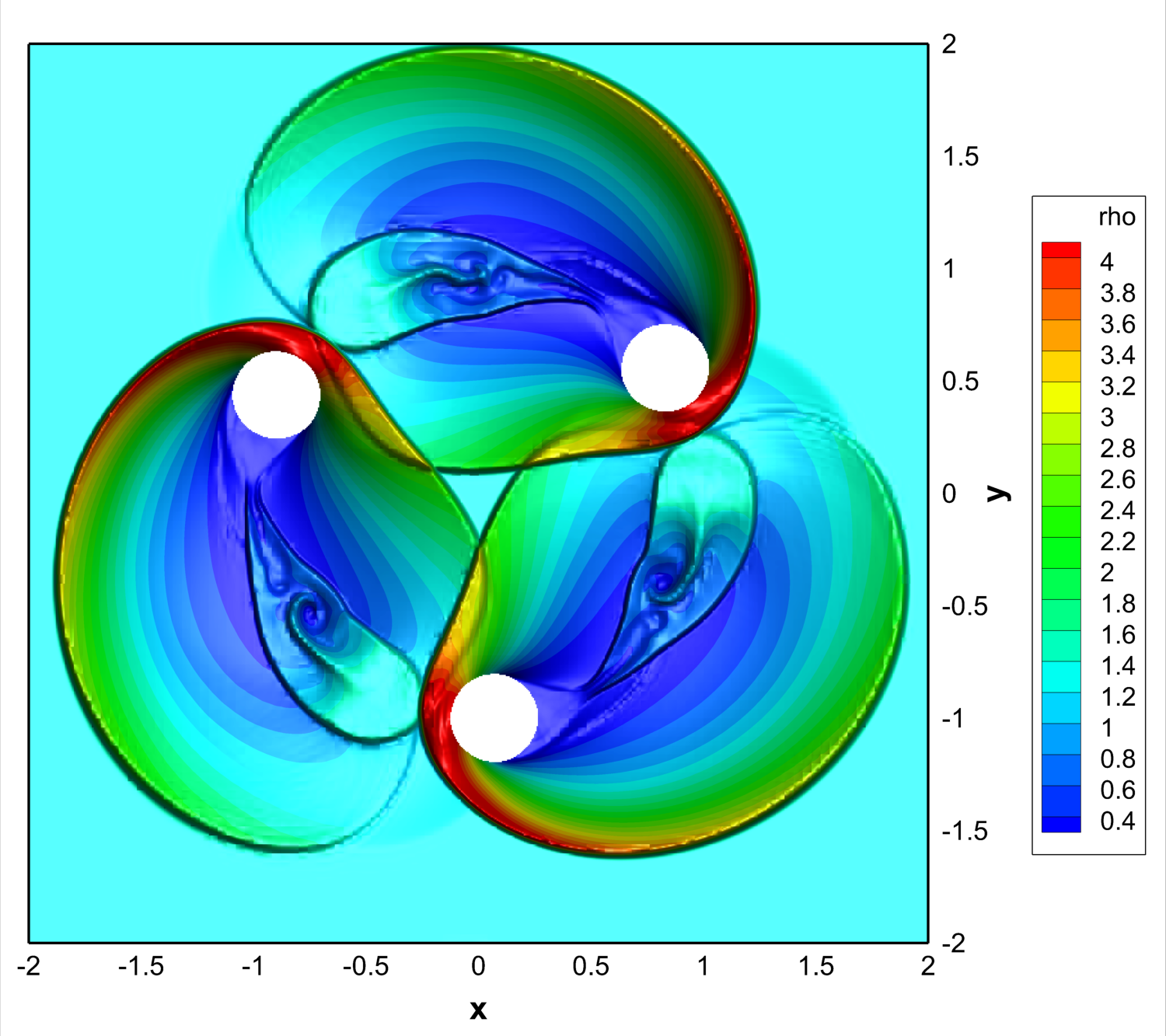}  \\ 
			\includegraphics[width=0.48\textwidth]{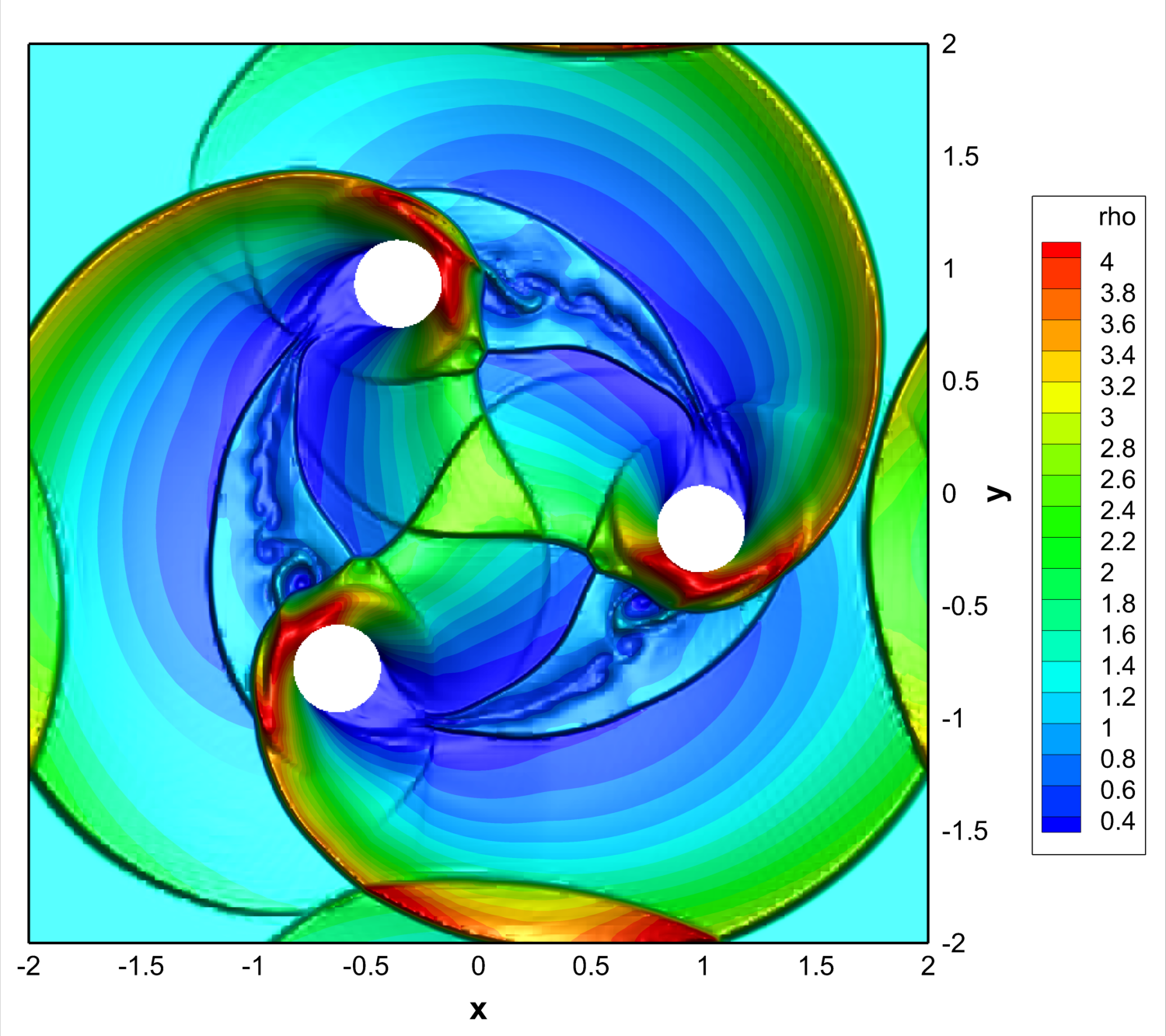}  & 
			\includegraphics[width=0.48\textwidth]{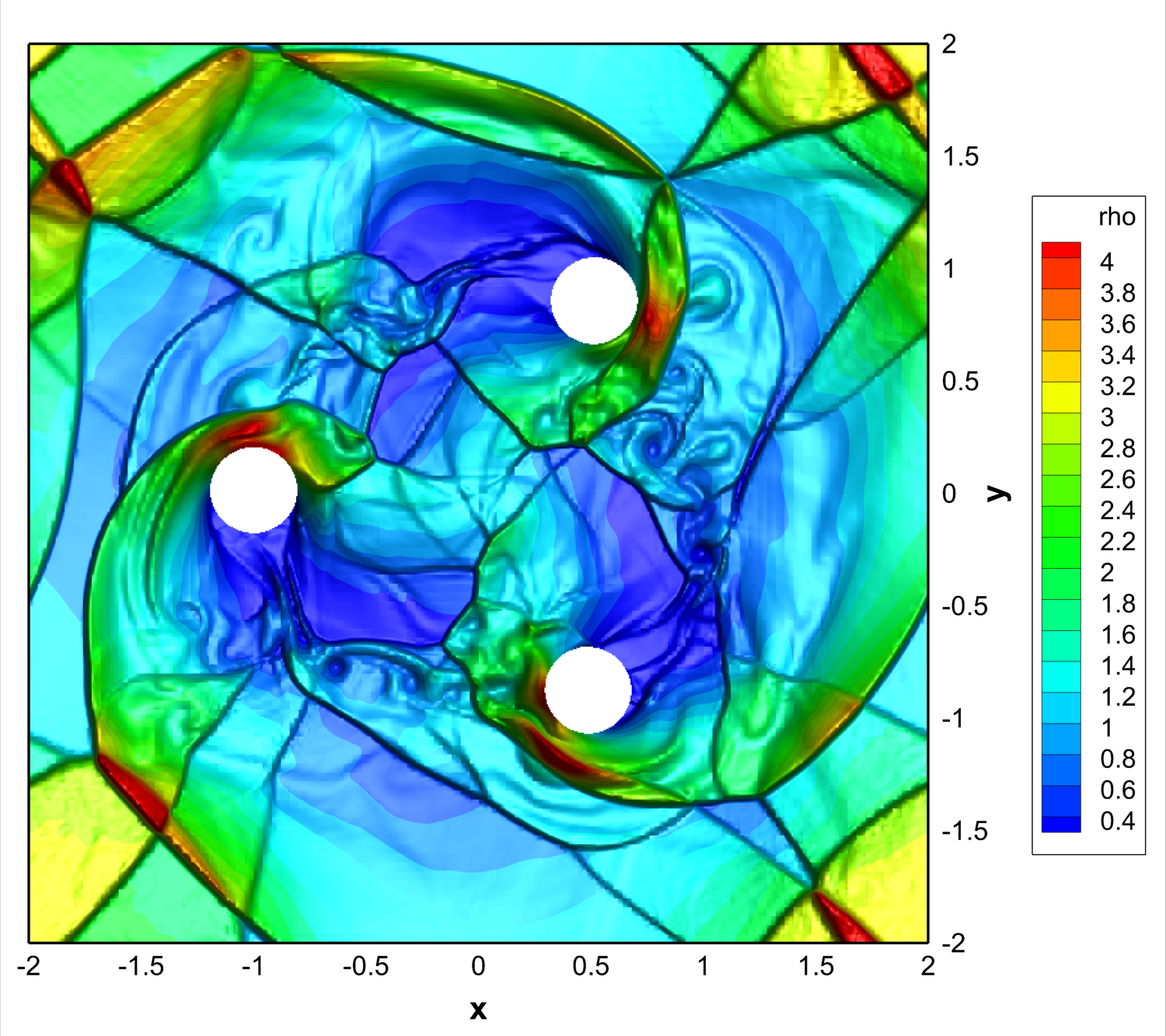}    
		\end{tabular}
	\end{center}
	\caption{Three cylinders rotating at supersonic speed in a compressible gas. Output times from top left to bottom right: $t=0.35$, $t=0.5$, $t=0.75$ and $t=1.75$.} 
	\label{fig.cyl}
\end{figure}

\begin{figure}
	\begin{center} 
		\begin{tabular}{cc}
			\includegraphics[width=0.48\textwidth]{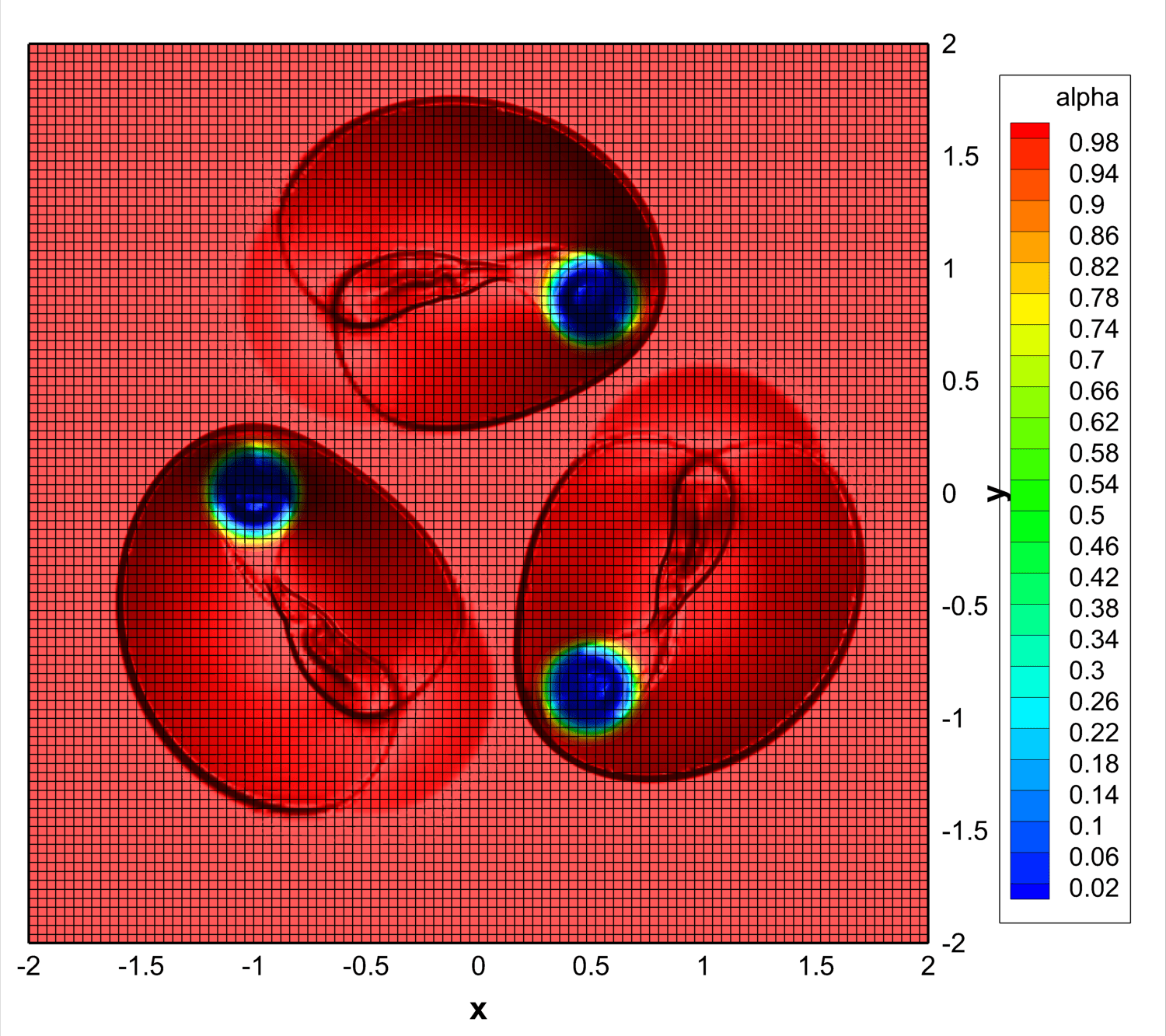}  & 
			\includegraphics[width=0.48\textwidth]{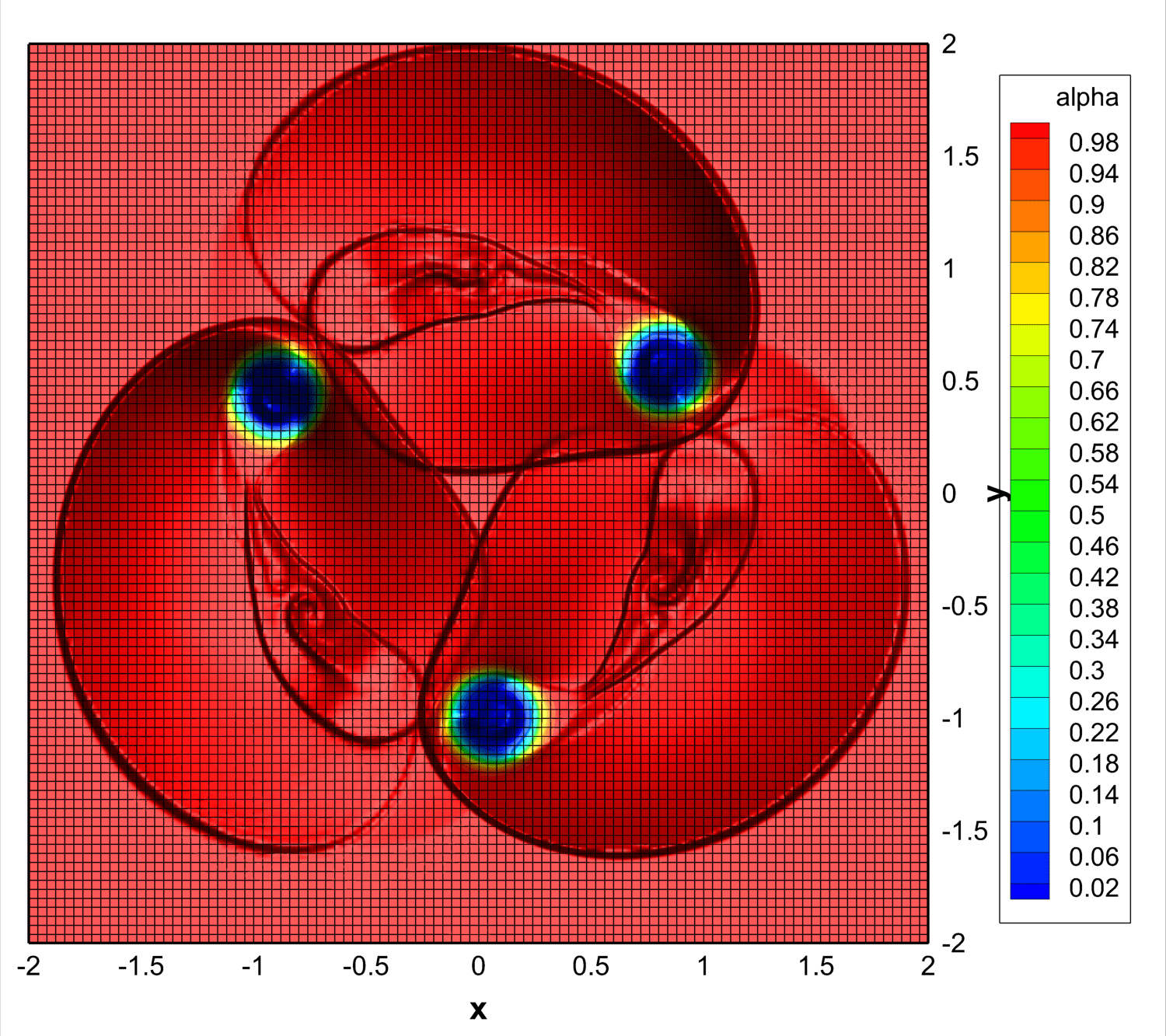}  \\ 
			\includegraphics[width=0.48\textwidth]{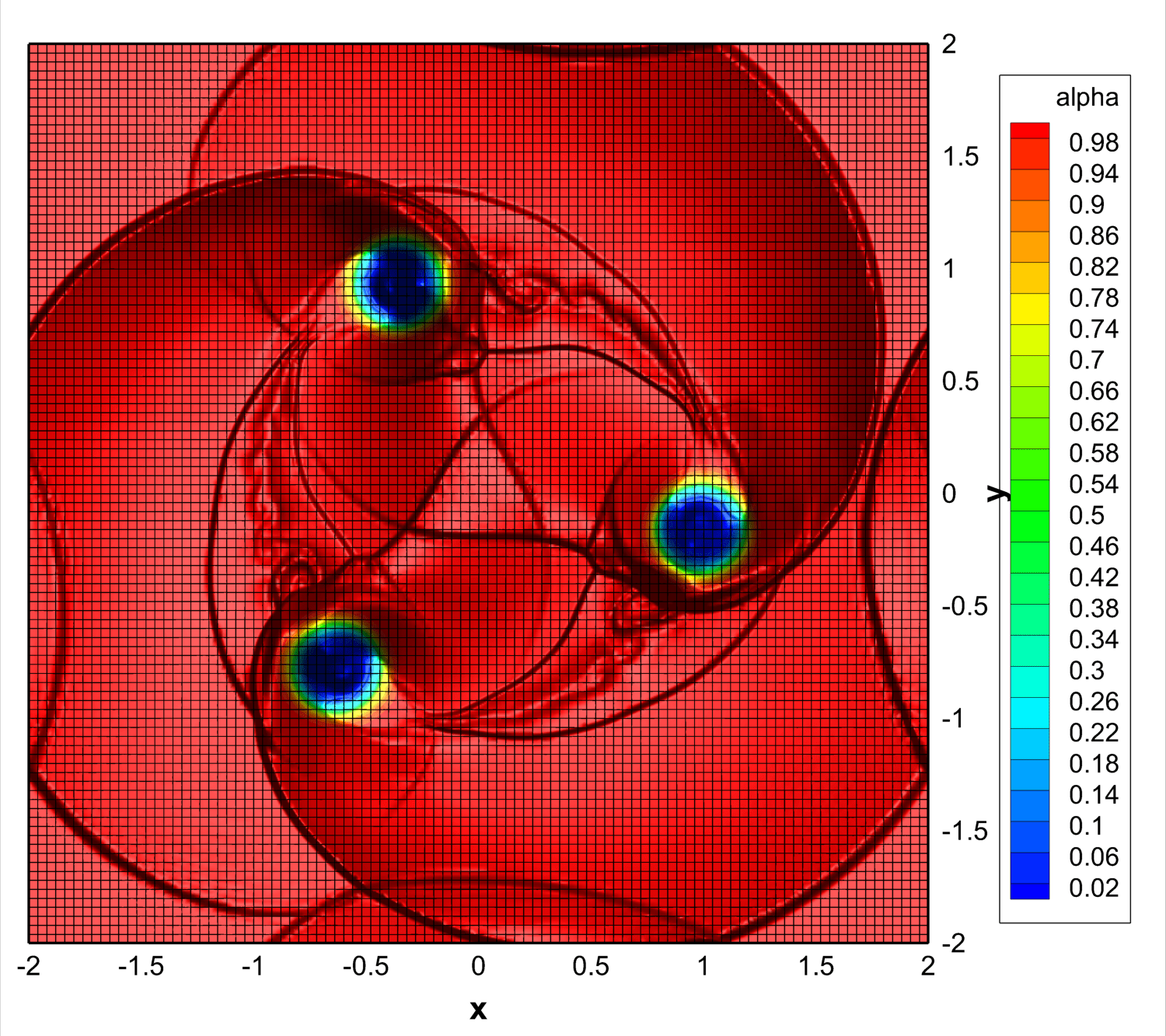}  & 
			\includegraphics[width=0.48\textwidth]{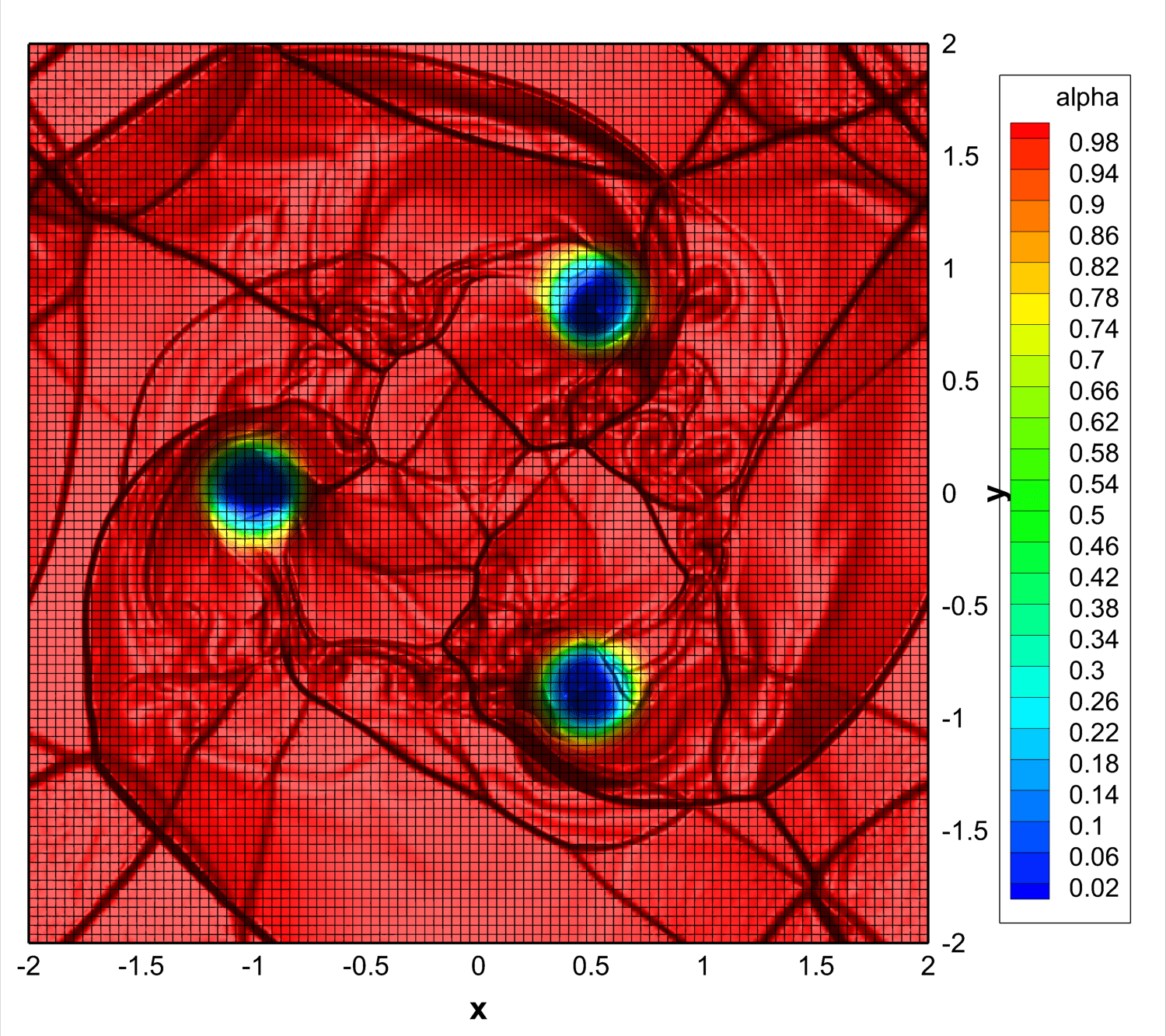}    
		\end{tabular}
	\end{center}
	\caption{Volume fraction function $\alpha$ for the three cylinders rotating at supersonic speed in a compressible gas. Output times from top left to bottom right: $t=0.35$, $t=0.5$, $t=0.75$ and $t=1.75$.} 
	\label{fig.cyl.alpha}
\end{figure}

This last example is only meant to be a showcase to illustrate the applicability of the diffuse interface model proposed in this paper. We choose the computational domain  $\Omega = [-2,+2] \times [-2, +2]$, which is discretized with a Cartesian grid of $100 \times 100$ equidistant ADER-DG elements. The polynomial degree of the basis functions is chosen as $N=3$. Periodic boundary conditions are applied everywhere. The initial condition for the fluid is $\rho = 1.4$, $\mathbf{v}=0$ and $p=1$, while the solid velocity is initialized as a rigid body rotation with  
$\mathbf{v}_s = \boldsymbol{\omega} \times \mathbf{x}$, with $\boldsymbol{\omega}=(0,0,-3)$. We consider three circular solid bodies of radius $R_c=0.2$, whose centers $\mathbf{x}_i$ are initially located on a circle of radius $R_0$, i.e. $\mathbf{x}_i = R_0 \left( \cos(\varphi_i), \sin(\varphi_i) \right)$ with 
$\varphi_i = \frac{i}{3} \cdot 2 \pi$. With this configuration, the centers of the circles move at a supersonic Mach number of $M=3$. The entire test problem is setup by simply defining the solid velocity field and by setting inside each circular solid body $\alpha = \epsilon$, while in the rest of the domain, we choose 
$\alpha = 1-\epsilon$, with $\epsilon = 10^{-3}$. Nothing else needs to be done. Simulations are carried out until $t=2.0$. The density contour colors at different output times are depicted in Figure \ref{fig.cyl}, where we have blanked regions where $\alpha<0.5$. One can clearly see the wake generated behind each cylinder, as well as the bow shock in front of each cylinder. At later times, the wakes and the bow shocks of two consecutive cylinders interact with each other. \textcolor{black}{In order to assess the level of numerical diffusion present in the volume fraction function $\alpha$, we also show the temporal evolution of $\alpha$ in Figure \ref{fig.cyl.alpha}. It can be seen that the interface is resolved within 3-4 DG cells, which is much more than in the case of the stationary blunt body shown in the previous section. In this context, less dissipative Riemann solvers, such as those forwarded in \cite{OsherNC,HLLEMNC}, and better slope limiters \cite{toro-book} need to be used in the future, also in combination with AMR \cite{AMR3DCL,DGLimiter2}. }  

One last time we stress that the entire setup of this test problem is very simple and that this very complex flow field and the description of the moving solid bodies is directly done inside the PDE system and only by prescribing a simple scalar volume fraction function. No moving body-fitted mesh needs to be generated.

\section{Conclusion} 
\label{sec.conclusion}

In this paper we have introduced a new and very simple diffuse interface model for the simulation of compressible flows around fixed and moving rigid bodies of arbitrary shape. The proposed approach is so simple that all simulations can be carried out on uniform Cartesian grids. The geometry of each solid body is simply defined by a scalar volume fraction function, which assumes a value close to zero inside the solids and close to unity inside the compressible fluid. The computational setup can be done in a fully automatic manner, without the need to generate a body-fitted structured or unstructured mesh, which can become very time consuming for complex geometries. The work presented in this paper is a natural extension of the simple diffuse interface approach introduced in \cite{DIM2D,DIM3D,Gaburro2018,Tavelli2019} and \cite{FrontierADERGPR}. The presented numerical algorithm with the underlying diffuse interface model is a consequent application of Godunov's shock capturing ideas to the context of moving material interfaces. Rather than tracking the material interface via a moving mesh, one simply needs to evolve in time the additional scalar quantity $\alpha$, which contains all the necessary information about the geometry of the problem.  

In the first part of this paper we have proven, via Riemann invariants and via generalized Rankine-Hugoniot conditions derived from the DLM theory \cite{DLMtheory}, that the normal component of the fluid velocity at the material interface assumes automatically the value of the normal component of the solid velocity when the fluid volume fraction $\alpha$ jumps from unity to zero, i.e. the \textit{non-penetration boundary condition} is naturally satisfied by the model. 
The proof based on generalized Rankine-Hugoniot conditions supposes a simple straight line segment path, which  therefore also justifies its use inside path-conservative schemes, as those employed in this paper. 

Future developments will concern the generalization to fluids interacting with elastic solids, via a combination of the compressible multi-phase model of Romenski et al. \cite{Rom1998,RomenskiTwoPhase2007,RomenskiTwoPhase2010} and the equations of hyperelasticity in Euler coordinates of Godunov, Peshkov and Romenski (GPR model), forwarded in \cite{GodunovRomenski72,GodRom2003,PeshRom2014,GPR1,GPR2}. We also plan to include surface tension effects via the new hyperbolic surface tension model of Gavrilyuk et al., see \cite{HypSurfTension,SHTCSurfaceTension}.  

\section*{Acknowledgments}
This research was funded by the European Union's Horizon 2020 Research and Innovation Programme under the project \textit{ExaHyPE}, grant no. 671698 (call FETHPC-1-2014). 

The authors would like to acknowledge PRACE for awarding access to the SuperMUC 
supercomputer based in Munich, Germany at the Leibniz Rechenzentrum (LRZ). 

M.D. acknowledges the financial support received from the Italian Ministry of Education, University and Research (MIUR) in the frame of the Departments of Excellence Initiative 2018--2022 attributed to DICAM of the University of Trento (grant L. 232/2016) and in the frame of the PRIN 2017 project. M.D. has also received funding from the University of Trento via the  \textit{Strategic Initiative Modeling and Simulation}. 

E.G. has been financed by a national mobility grant for young researchers in Italy, funded by GNCS-INdAM and acknowledges the support given by the University of Trento through the \textit{UniTN Starting Grant} initiative.


\noindent \section*{In memoriam}

\noindent This paper is dedicated to the memory of Dr. Douglas Nelson Woods ($^*$January 11\textsuperscript{th} 1985 - $\dagger$September 11\textsuperscript{th} 2019),
promising young scientist and post-doctoral research fellow at Los Alamos National Laboratory. 
Our thoughts and wishes go to his wife Jessica, to his parents Susan and Tom, to his sister Rebecca and to his brother Chris, whom he left behind.

\clearpage 

\bibliographystyle{plain}
\bibliography{references}


\end{document}